\documentclass[12pt]{article}
\usepackage{amsfonts}
\usepackage{amssymb}
\usepackage{amsmath}
\usepackage[numbers,sort&compress]{natbib}
\usepackage[colorlinks,citecolor=blue,urlcolor=blue]{hyperref}
\usepackage{graphicx,indentfirst,tabularx}
\usepackage[T1]{fontenc}
\usepackage{enumitem}
\usepackage{mathrsfs}

\newtheorem{theorem}{Theorem}

\newtheorem{corollary}[theorem]{Corollary}

\newtheorem{definition}[theorem]{Definition}
\newtheorem{example}[theorem]{Example}

\newtheorem{lemma}[theorem]{Lemma}

\newtheorem{proposition}[theorem]{Proposition}
\newtheorem{remark}[theorem]{Remark}

\newenvironment{proof}[1][Proof]{\noindent\textbf{#1.} }{\ \rule{0.5em}{0.5em}}
\input{tcilatex}
\begin{document}

\title{Stochastic Multivariate Taylor Measure Functions}
\author{Athanasios C. Micheas \thanks{%
Department of Statistics, University of Missouri, 146 Middlebush Hall,
Columbia, MO 65211-6100, USA, email: micheasa@missouri.edu}}
\maketitle

\begin{abstract}
We introduce stochastic multivariate Taylor measure functions (SMTMFs),
which extend the multivariate Taylor measure function framework through
random Taylor coefficients and a possibly random Taylor parameter. We
establish their measurability properties and derive sufficient conditions
for the associated processes to be martingales, including criteria based on
conditional covariance and continuous semimartingale decompositions.

We construct standard Brownian motion as an almost surely uniformly
convergent limit of finite SMTMFs using the L\'{e}vy-Ciesielski
representation. We also show that the completion of the finite stochastic
MTMF space under the expected uniform norm contains Brownian motion,
Brownian bridges, fractional Brownian motion with Hurst parameter $H\in(0,1)$%
, and Ornstein-Uhlenbeck processes with integrable initial conditions.

Our main result characterizes this completion as the Bochner space 
\begin{equation*}
L^1(\Omega;C([0,T])),
\end{equation*}
under an admissibility condition on jointly measurable stochastic Taylor
coefficients. The proof uses dyadic piecewise linear interpolation to
approximate integrable continuous stochastic processes in the expected
uniform norm. Thus, the completed SMTMF space encompasses the full class of
integrable continuous random functions on a compact time interval.
\end{abstract}



\section{Introduction}

Taylor measure functions provide a function space framework based on Taylor
measure representations. The theory of Taylor measures and its multivariate
function space extension were developed in \cite{micheas2025taylor} and \cite%
{micheas2026MTMF}, respectively. These constructions offer a way to
represent functions through Taylor coefficients and a Taylor parameter,
while retaining the structure of the underlying Taylor measure. The purpose
of this paper is to extend this framework to random functions by allowing
both the Taylor coefficients and the Taylor parameter to be stochastic. This
leads to the class of stochastic multivariate Taylor measure functions
(SMTMFs), which connects stochastic Taylor measures with multivariate Taylor
measure functions.

More precisely, let $(\Omega ,\mathcal{A},\mathbb{P})$ be a probability
space, let $p\in \mathbb{N}^{+}=\{1,2,...\}$, and let $\mathbb{M}\in \{%
\mathbb{R},\mathbb{C}\}$. We consider jointly measurable functions 
\begin{equation*}
g:\Omega \times \mathbb{M}^{p}\longrightarrow \mathbb{M},
\end{equation*}%
and 
\begin{equation*}
a_{n}:\Omega \times \mathbb{M}^{p}\longrightarrow \mathbb{M},\qquad n\in 
\mathbb{N}_{0}=\{0,1,2,...\},
\end{equation*}%
which are called the Taylor parameter and Taylor coefficients, respectively.
For a set $B\subseteq \mathbb{N}_{0}$, the corresponding stochastic
multivariate Taylor measure function is given by 
\begin{equation*}
F_{B}(\omega ,x)=\sum_{n\in B}a_{n}(\omega ,x)\frac{g(\omega ,x)^{n}}{n!}%
,\qquad x\in \mathbb{M}^{p},
\end{equation*}%
whenever the associated Taylor measure has finite total variation. Thus, the
defining representation retains the structure of a Taylor measure while
allowing dependence on both the random outcome and the spatial variable.

The first part of the paper develops the fundamental measurability and
probabilistic properties of this construction. In particular, we establish
the existence of jointly measurable versions of SMTMFs under the stated
assumptions. An important measure theoretic point is the distinction between
almost sure convergence at each fixed spatial input and simultaneous almost
sure convergence over an uncountable spatial domain. The former does not, by
itself, imply the latter. We therefore formulate the measurability results
in terms of jointly measurable versions and identify additional conditions
under which uniform convergence on compact sets can be obtained.

We next investigate the martingale properties of SMTMFs. Since the defining
Taylor coefficients and Taylor parameter may both be stochastic and time
dependent, the martingale property of the resulting function cannot in
general be deduced merely from the martingale property of the individual
coefficients. Rather, it is the Taylor terms 
\begin{equation*}
Y_{n,t}(x)=a_{n}(t,x)\frac{g(t,x)^{n}}{n!}
\end{equation*}%
that play a central role. We establish sufficient conditions under which a
summation of integrable martingale Taylor terms is itself a martingale. We
also derive a conditional covariance criterion for the product of a
coefficient process and a power of the Taylor parameter. In the continuous
semimartingale setting, It\^{o}'s formula yields an explicit drift condition
involving the drift of the coefficient, the drift of the Taylor parameter,
their quadratic covariation, and the quadratic variation of the Taylor
parameter. This condition describes how the different stochastic
contributions can compensate one another so that an individual Taylor term
is a local martingale, and, under an additional class $D$ condition, a true
martingale.

A further objective is to demonstrate that the SMTMF construction can
represent fundamental stochastic processes. We begin with an explicit
construction of standard Brownian motion based on the L\'{e}vy--Ciesielski
representation. This representation expresses Brownian motion through
independent standard Gaussian random variables and deterministic Schauder
functions. We show that its finite approximations admit finite stochastic
MTMF representations and converge almost surely uniformly to a process with
continuous sample paths, centered Gaussian finite dimensional distributions,
and covariance function $\min (s,t)$. The limiting process is therefore
standard Brownian motion. This construction provides an explicit example of
how a classical stochastic process arises as a uniform limit of finite
elements of the proposed function space.

We then introduce a completion of the finite stochastic MTMF space under an
appropriate expected uniform norm. This completion makes it possible to
treat limits of finite stochastic Taylor representations as elements of a
well defined function space. In particular, we establish the inclusion of
Brownian motion, Brownian bridges, fractional Brownian motion with Hurst
parameter $H\in (0,1)$, and Ornstein--Uhlenbeck processes with integrable
initial conditions. The arguments use the continuity of the sample paths and
integrability of the uniform norm, together with finite stochastic MTMF
approximations constructed by piecewise linear interpolation. These results
show that the completed space accommodates important continuous stochastic
processes arising in different probabilistic settings.

The principal result of the paper is a characterization of the completed
stochastic MTMF space. Let $C([0,T])$ denote the Banach space of real-valued
continuous functions on a compact interval $[0,T]$, equipped with the
supremum norm. We equip the finite stochastic MTMF space with the expected
uniform norm 
\begin{equation*}
\Vert F\Vert _{1,\infty }=\left[ \sup_{0\leq t\leq T}|F(\omega ,t)|\right] .
\end{equation*}%
Under the admissibility conditions specified in the paper, in particular the
allowance of jointly measurable stochastic Taylor coefficients depending on
both $\omega $ and $t$, we prove that its completion satisfies 
\begin{equation*}
\overline{\mathscr S}_{1}=L^{1}(\Omega ;C([0,T])).
\end{equation*}%
The key step in the proof is the construction, for an arbitrary integrable
continuous stochastic process, of finite stochastic MTMF approximations by
dyadic piecewise linear interpolation. The approximations converge almost
surely uniformly and also in the expected uniform norm. Conversely, every
finite stochastic MTMF belongs to the corresponding Bochner space, so its
completion is contained in that space.

This characterization is the central contribution of the paper. It
establishes that the completed stochastic MTMF space is not merely a class
containing selected Gaussian processes or particular stochastic
constructions. Under the specified coefficient admissibility convention, it
is precisely the space of integrable continuous random functions on the
compact time interval. It therefore provides a general functional analytic
framework in which continuous stochastic processes can be represented and
approximated through finite stochastic Taylor measure functions. We
emphasize that this identification concerns the completion of the finite
stochastic MTMF space; it does not assert that every element of the Bochner
space has a finite stochastic MTMF representation. The admissibility of
jointly measurable random coefficients is essential to the interpolation
argument.

The remainder of the paper is organized as follows. Section 2 introduces
stochastic multivariate Taylor measure functions and develops their
fundamental measurability and martingale properties, including conditional
covariance and semimartingale criteria. Section 3 gives an explicit
construction of standard Brownian motion through finite SMTMFs using the L%
\'{e}vy--Ciesielski representation. Section 4 introduces the completion of
the stochastic MTMF space and establishes the inclusion of Brownian motion,
Brownian bridges, fractional Brownian motion, and Ornstein--Uhlenbeck
processes. Section 5 establishes the characterization of the completed space
in terms of the Bochner space $L^1(\Omega;C([0,T]))$. Finally, Section 6
presents concluding remarks and discusses directions for further research.

\section{Definitions and Properties}

We first recall the concepts of finite Taylor measures and stochastic Taylor
measures that were introduced\ and studied in \cite{micheas2025taylor}.

\begin{definition}[Finite Taylor Measure]
Let $\mathbb{M}\in \{\mathbb{R},\mathbb{C}\}$ and $\mathbb{N}%
_{0}=\{0,1,2,\ldots \}$. For $\gamma \in \mathbb{M}$ and $%
a=(a_{0},a_{1},\ldots )\in \mathbb{M}^{\mathbb{N}_{0}}$, define 
\begin{equation*}
T_{\gamma ,a}(B)=\sum_{n\in B}a_{n}\frac{\gamma ^{n}}{n!},\qquad B\in 
\mathcal{B}(\mathbb{N}_{0}).
\end{equation*}%
The collection of all such set functions satisfying 
\begin{equation*}
\sum_{n=0}^{\infty }|a_{n}|\frac{|\gamma |^{n}}{n!}<\infty
\end{equation*}%
is denoted by $\mathcal{T}^{F}$.
\end{definition}

\begin{definition}[Stochastic Taylor Measure]
Let $(\Omega ,\mathcal{A},P)$ be a probability space, and 
\begin{equation*}
a:\Omega \rightarrow \mathbb{M}^{\mathbb{N}_{0}},\qquad \gamma :\Omega
\rightarrow \mathbb{M}
\end{equation*}%
be measurable mappings. If 
\begin{equation*}
\sum_{n=0}^{\infty }|a_{n}(\omega )|\frac{|\gamma (\omega )|^{n}}{n!}<\infty
\quad \text{almost surely},
\end{equation*}%
then the random set function 
\begin{equation*}
X(\omega )(B)=\sum_{n\in B}a_{n}(\omega )\frac{\gamma (\omega )^{n}}{n!}
\end{equation*}%
is called a Stochastic Taylor Measure.
\end{definition}

On the other hand, the multivariate Taylor measure function space (MTMF) was
introduced and studied in \cite{micheas2026MTMF}, which is a function space
with purely deterministic elements.

\begin{definition}[Multivariate Taylor Measure Function]
Let $p\in \mathbb{N}^{+}=\left\{ 1,2,...\right\} $ and $g,a_{n}:\mathbb{M}%
^{p}\rightarrow \mathbb{M}$ be measurable functions. For $B\in \mathcal{B}(%
\mathbb{N}_{0})$, define 
\begin{equation*}
f^{B}(x)=\sum_{n\in B}a_{n}(x)\frac{g(x)^{n}}{n!},\qquad x\in \mathbb{M}^{p},
\end{equation*}%
provided that 
\begin{equation*}
\sum_{n=0}^{\infty }|a_{n}(x)|\frac{|g(x)|^{n}}{n!}<\infty .
\end{equation*}%
Then $f^{B}$ is called a Multivariate Taylor Measure Function.
\end{definition}

Based on the above definitions, we are now ready to propose a general
stochastic version.

\begin{definition}[Stochastic Multivariate Taylor Measure Function]
Let $(\Omega ,\mathcal{A},\mathbb{P})$ be a probability space, $p\in \mathbb{%
N}^{+}$, and $\mathbb{M}\in \{\mathbb{R},\mathbb{C}\}$.

Consider jointly measurable functions 
\begin{equation*}
g:\Omega \times \mathbb{M}^{p}\rightarrow \mathbb{M}
\end{equation*}%
which will be called the Taylor parameter, and 
\begin{equation*}
a_{n}:\Omega \times \mathbb{M}^{p}\rightarrow \mathbb{M},\qquad n\in \mathbb{%
N}_{0},
\end{equation*}%
which will be called the Taylor coefficients. For every $B\in \mathcal{B}(%
\mathbb{N}_{0})$, define 
\begin{equation*}
F^{B}(\omega ,x)=\sum_{n\in B}a_{n}(\omega ,x)\frac{g(\omega ,x)^{n}}{n!}.
\end{equation*}%
Suppose that 
\begin{equation*}
\sum_{n=0}^{\infty }|a_{n}(\omega ,x)|\frac{|g(\omega ,x)|^{n}}{n!}<\infty
\end{equation*}%
for almost every $\omega $, for each fixed $x\in \mathbb{M}^{p}$.

Then $F^{B}$ is called a Stochastic Multivariate Taylor Measure Function
(SMTMF).
\end{definition}

The collection of all SMTMFs is defined next.

\begin{definition}[General Space of Stochastic Taylor Measure Functions]
The space of all SMTMFs satisfying the preceding definition is denoted by 
\begin{equation*}
\mathscr{S}_{p}(\mathbb{M})=\left\{ F^{B}:F^{B}(\omega ,x)=\sum_{n\in
B}a_{n}(\omega ,x)\frac{g(\omega ,x)^{n}}{n!}\right\} ,
\end{equation*}%
where the defining functions are jointly measurable and the corresponding
Taylor measures have finite total variation.
\end{definition}

\subsection{Joint Measurability of SMTMFs}

Let $(\Omega ,\mathcal{A},\mathbb{P})$ be a probability space, let $p\in 
\mathbb{N}^{+}$, and let $\mathbb{M}\in \{\mathbb{R},\mathbb{C}\}$. Consider
the measurable spaces 
\begin{equation*}
(\Omega ,\mathcal{A})\quad \text{and}\quad (\mathbb{M}^{p},\mathcal{B}(%
\mathbb{M}^{p})),
\end{equation*}%
and suppose that 
\begin{equation*}
g:\Omega \times \mathbb{M}^{p}\longrightarrow \mathbb{M},
\end{equation*}%
and 
\begin{equation*}
a_{n}:\Omega \times \mathbb{M}^{p}\longrightarrow \mathbb{M},\qquad n\in 
\mathbb{N}_{0},
\end{equation*}%
are jointly measurable with respect to 
\begin{equation*}
\mathcal{A}\otimes \mathcal{B}(\mathbb{M}^{p}),
\end{equation*}%
and $\mathcal{B}(\mathbb{M})$. For every $B\in \mathcal{B}(\mathbb{N}_{0})$,
define 
\begin{equation*}
F^{B}(\omega ,x)=\sum_{n\in B}a_{n}(\omega ,x)\frac{g(\omega ,x)^{n}}{n!}.
\end{equation*}%
Assume that, for every fixed $x\in \mathbb{M}^{p}$, such that 
\begin{equation*}
\sum_{n=0}^{\infty }|a_{n}(\omega ,x)|\frac{|g(\omega ,x)|^{n}}{n!}<\infty
\quad \text{almost surely}.
\end{equation*}%
We establish the existence of a jointly measurable version of the stochastic
Taylor measure function in the following.

\begin{proposition}[Joint Measurability]
Every element $F^{B}\in \mathscr{S}_{p}(\mathbb{M})$ admits a jointly
measurable version on $\Omega \times \mathbb{M}^{p}$. Furthermore, for every
fixed $x\in \mathbb{M}^{p}$, the mapping 
\begin{equation*}
F^{B}(\cdot ,x):\Omega \longrightarrow \mathbb{M},
\end{equation*}%
is a random variable.
\end{proposition}

\begin{proof}
We divide the proof into seven steps.

\noindent \textbf{Step 1. Definition of the partial sums}

For every $N\in \mathbb{N}_{0}$, define 
\begin{equation*}
F_{N}^{B}(\omega ,x)=\sum_{\substack{ n\in B  \\ n\leq N}}a_{n}(\omega ,x)%
\frac{g(\omega ,x)^{n}}{n!}.
\end{equation*}%
Since $B\subseteq \mathbb{N}_{0}$, the sum contains only finitely many terms
for each $N$. We first establish that each $F_{N}^{B}$ is jointly measurable.

\noindent \textbf{Step 2. Measurability of each summand}

By assumption, the mappings 
\begin{equation*}
a_{n}:\Omega \times \mathbb{M}^{p}\rightarrow \mathbb{M},
\end{equation*}%
and 
\begin{equation*}
g:\Omega \times \mathbb{M}^{p}\rightarrow \mathbb{M},
\end{equation*}%
are jointly measurable. For each fixed $n\in \mathbb{N}_{0}$, consider 
\begin{equation*}
h_{n}(\omega ,x)=a_{n}(\omega ,x)g(\omega ,x)^{n}.
\end{equation*}%
The mapping 
\begin{equation*}
\varphi _{n}:\mathbb{M}\rightarrow \mathbb{M},\qquad \varphi _{n}(z)=z^{n},
\end{equation*}%
is continuous. Therefore, the composition 
\begin{equation*}
(\omega ,x)\longmapsto g(\omega ,x)^{n},
\end{equation*}%
is measurable. Furthermore, the multiplication mapping 
\begin{equation*}
m:\mathbb{M}\times \mathbb{M}\rightarrow \mathbb{M},\qquad m(u,v)=uv,
\end{equation*}%
is continuous, and it follows that 
\begin{equation*}
(\omega ,x)\longmapsto a_{n}(\omega ,x)g(\omega ,x)^{n},
\end{equation*}%
is measurable. Since division by the positive constant $n!$ preserves
measurability, we conclude that 
\begin{equation*}
(\omega ,x)\longmapsto a_{n}(\omega ,x)\frac{g(\omega ,x)^{n}}{n!},
\end{equation*}%
is jointly measurable for every $n\in \mathbb{N}_{0}$. For $n=0$, we use the
usual convention $z^{0}=1$, including when $z=0$.

\noindent \textbf{Step 3. Measurability of the partial sums}

Recall that 
\begin{equation*}
F_{N}^{B}(\omega ,x)=\sum_{\substack{ n\in B  \\ n\leq N}}a_{n}(\omega ,x)%
\frac{g(\omega ,x)^{n}}{n!}.
\end{equation*}%
By Step 2, every summand is jointly measurable. Since a finite sum of
measurable functions is measurable, we obtain 
\begin{equation*}
F_{N}^{B}:\Omega \times \mathbb{M}^{p}\rightarrow \mathbb{M},
\end{equation*}%
jointly measurable for every $N\in \mathbb{N}_{0}$.

\noindent \textbf{Step 4. Establishment of the convergence set}

Define 
\begin{equation*}
D=\left\{ (\omega ,x)\in \Omega \times \mathbb{M}^{p}:\lim_{N\rightarrow
\infty }F_{N}^{B}(\omega ,x)\text{ exists in }\mathbb{M}\right\} .
\end{equation*}%
We show that $D$ is measurable. Since $\mathbb{M}$ is either $\mathbb{R}$ or 
$\mathbb{C}$, it is a complete metric space under its usual metric.
Consequently, a sequence in $\mathbb{M}$ converges if and only if it is a
Cauchy sequence. Thus, we can write 
\begin{equation*}
D=\bigcap_{k=1}^{\infty }\bigcup_{N=0}^{\infty }\bigcap_{r,s\geq N}\left\{
(\omega ,x):|F_{r}^{B}(\omega ,x)-F_{s}^{B}(\omega ,x)|<\frac{1}{k}\right\} .
\end{equation*}%
By Step 3, the functions $F_{r}^{B}$ and $F_{s}^{B}$ are jointly measurable.
Hence, it holds that%
\begin{equation*}
(\omega ,x)\longmapsto |F_{r}^{B}(\omega ,x)-F_{s}^{B}(\omega ,x)|,
\end{equation*}%
is jointly measurable. Therefore, every set appearing on the right-hand side
belongs to $\mathcal{A}\otimes \mathcal{B}(\mathbb{M}^{p})$. Since $\sigma $%
-algebras are closed under countable unions and intersections, we have 
\begin{equation*}
D\in \mathcal{A}\otimes \mathcal{B}(\mathbb{M}^{p}).
\end{equation*}

\noindent \textbf{Step 5. Construction of a jointly measurable version}

Define 
\begin{equation*}
\widetilde{F}^{B}(\omega ,x)=%
\begin{cases}
\lim_{N\rightarrow \infty }F_{N}^{B}(\omega ,x), & (\omega ,x)\in D, \\[8pt] 
0, & (\omega ,x)\notin D.%
\end{cases}%
\end{equation*}%
For every $N\in \mathbb{N}_{0}$, define 
\begin{equation*}
G_{N}^{B}(\omega ,x)=\mathbf{1}_{D}(\omega ,x)F_{N}^{B}(\omega ,x),
\end{equation*}%
where $\mathbf{1}_{D}$ denotes the indicator function of $D$.

Since $D$ is measurable and $F_{N}^{B}$ is jointly measurable, $G_{N}^{B}$
is jointly measurable. Moreover, we have%
\begin{equation*}
\widetilde{F}^{B}(\omega ,x)=\lim_{N\rightarrow \infty }G_{N}^{B}(\omega ,x).
\end{equation*}%
Indeed, if $(\omega ,x)\in D$, the sequence converges to the defining limit.
If $(\omega ,x)\notin D$, every $G_{N}^{B}(\omega ,x)$ equals zero. By the
pointwise limit theorem for measurable functions, $\widetilde{F}^{B}$ is
jointly measurable.

\noindent \textbf{Step 6. Identification with the original series almost
surely}

By the assumed absolute convergence condition, for every fixed $x\in \mathbb{%
M}^{p}$, we have%
\begin{equation*}
P\left( \left\{ \omega :\sum_{n=0}^{\infty }|a_{n}(\omega ,x)|\frac{%
|g(\omega ,x)|^{n}}{n!}<\infty \right\} \right) =1.
\end{equation*}%
Absolute convergence implies convergence of the series over every subset $%
B\subseteq \mathbb{N}_{0}$. Consequently, 
\begin{equation*}
P\left( \left\{ \omega :(\omega ,x)\in D\right\} \right) =1,
\end{equation*}%
and thus, for every fixed $x\in \mathbb{M}^{p}$, 
\begin{equation*}
\widetilde{F}^{B}(\omega ,x)=F^{B}(\omega ,x)\quad \text{almost surely}.
\end{equation*}%
Therefore, $\widetilde{F}^{B}$ is a jointly measurable version of the
original series.

\noindent \textbf{Step 7. Establishment of the random-variable property}

Fix an arbitrary $x\in \mathbb{M}^{p}$. By Step 5, the mapping 
\begin{equation*}
\omega \longmapsto \widetilde{F}^{B}(\omega ,x),
\end{equation*}%
is measurable with respect to $\mathcal{A}$ and $\mathcal{B}(\mathbb{M})$.
Since $(\Omega ,\mathcal{A},P)$ is a probability space, every measurable
mapping from $\Omega $ into $\mathbb{M}$ is an $\mathbb{M}$-valued random
variable. Hence, $\widetilde{F}^{B}(\cdot ,x)$ is a random variable for
every fixed $x\in \mathbb{M}^{p}$. Since $\widetilde{F}^{B}(\cdot
,x)=F^{B}(\cdot ,x)$ almost surely, the original series also defines a
random variable up to almost sure equivalence. This completes the proof.
\end{proof}

There is an important distinction between pointwise almost sure convergence
and simultaneous almost sure convergence. The assumption 
\begin{equation*}
\forall x\in \mathbb{M}^{p},\qquad P(\Omega _{x})=1
\end{equation*}%
does not necessarily imply 
\begin{equation*}
P\left( \bigcap_{x\in \mathbb{M}^{p}}\Omega _{x}\right) =1,
\end{equation*}%
because $\mathbb{M}^{p}$ is generally uncountable. Therefore, the
proposition establishes the existence of a jointly measurable version of the
stochastic Taylor measure function. It does not automatically establish that
the original series converges simultaneously for every $x$ on a single event
of probability one.

If simultaneous convergence is desired, a stronger assumption is needed. For
example, let $K\subseteq \mathbb{M}^{p}$ be compact and suppose that 
\begin{equation*}
P\left( \left\{ \omega :\sum_{n=0}^{\infty }\sup_{x\in K}\left\vert
a_{n}(\omega ,x)\frac{g(\omega ,x)^{n}}{n!}\right\vert <\infty \right\}
\right) =1.
\end{equation*}%
Then, by the Weierstrass M-test, the defining series converges absolutely
and uniformly on $K$, almost surely. If the individual summands are
continuous in $x$, the resulting sample paths are continuous on $K$ almost
surely.

\subsection{General Martingale Conditions for SMTMFs}

A STMTF can be a martingale under suitable conditions on its random Taylor
coefficients, its Taylor parameter, and the underlying filtration. To
establish the martingale property, we introduce a filtration that describes
the evolution of information over time.

Let $(\Omega ,\mathcal{A},\mathbb{P})$ be a probability space equipped with
a filtration $\mathbb{F}=(\mathcal{A}_{t})_{t\geq 0},$ where 
\begin{equation*}
\mathcal{A}_{s}\subseteq \mathcal{A}_{t}\quad \text{whenever }s\leq t.
\end{equation*}%
Let $p\in \mathbb{N}^{+}$, and consider an STMTF of the form 
\begin{equation*}
F_{t}^{B}(\omega ,x)=\sum_{n\in B}a_{n}(t,\omega ,x)\frac{g(t,\omega ,x)^{n}%
}{n!},
\end{equation*}%
where $B\in \mathcal{B}(\mathbb{N}_{0})$. For every fixed $x\in \mathbb{M}%
^{p}$, we obtain a stochastic process 
\begin{equation*}
F_{t}^{B}(x)=F_{t}^{B}(\cdot ,x).
\end{equation*}

Our objective is to determine conditions under which this process is a
martingale with respect to the filtration $\mathbb{F}$. For completeness,
recall the definition of a martingale

\begin{definition}[Martingale]
Let $(\Omega ,\mathcal{A},\mathbb{F},P)$ be a filtered probability space. A
stochastic process $(X_{t})_{t\geq 0}$ is called a martingale with respect
to $\mathbb{F}$ if the following three conditions hold:

\begin{enumerate}
\item \textbf{Adaptedness:} For every $t\geq 0$, $X_t$ is $\mathcal{A}_t$%
-measurable.

\item \textbf{Integrability:} For every $t\geq 0$, $\mathbb{E}%
[|X_{t}|]<\infty .$

\item \textbf{Conditional expectation property:} For every $0\leq s\leq t$, 
\begin{equation*}
\mathbb{E}[X_t\mid\mathcal{A}_s]=X_s \quad\text{almost surely}.
\end{equation*}
\end{enumerate}
\end{definition}

For complex-valued processes, conditional expectation is defined
componentwise through the real and imaginary parts. We first formulate a
general sufficient condition that connect STMTFs with martingales.

\begin{theorem}[Martingale Characterization]
Let 
\begin{equation*}
F_{t}^{B}(x)=\sum_{n\in B}a_{n}(t,x)\frac{g(t,x)^{n}}{n!},
\end{equation*}%
where the random coefficients and Taylor parameter are defined on the
filtered probability space $(\Omega ,\mathcal{A},\mathbb{F},P).$ Fix $x\in 
\mathbb{M}^{p}$, and suppose the following conditions hold:

\begin{enumerate}
\item For every $t\geq 0$ and $n\in \mathbb{N}_{0}$, the functions $%
a_{n}(t,x)$ and $g(t,x)$ are $\mathcal{A}_{t}$-measurable.

\item For every $t\geq 0$, we have%
\begin{equation*}
\mathbb{E}\left[ \left\vert \sum_{n\in B}a_{n}(t,x)\frac{g(t,x)^{n}}{n!}%
\right\vert \right] <\infty .
\end{equation*}

\item For every $0\leq s\leq t$, it holds%
\begin{equation*}
\mathbb{E}\left[ \sum_{n\in B}a_{n}(t,x)\frac{g(t,x)^{n}}{n!}\middle|%
\mathcal{A}_{s}\right] =\sum_{n\in B}a_{n}(s,x)\frac{g(s,x)^{n}}{n!}\quad 
\text{a.s.}
\end{equation*}
\end{enumerate}

Then $\left( F_{t}^{B}(x)\right) _{t\geq 0}$ is a martingale with respect to 
$\mathbb{F}$.
\end{theorem}

\begin{proof}
We verify the three defining properties of a martingale.

\noindent \textbf{Step 1. Adaptedness.}

By assumption, for every $t\geq 0$ and $n\in \mathbb{N}_{0}$, the random
variables 
\begin{equation*}
a_{n}(t,x),\qquad g(t,x),
\end{equation*}%
are $\mathcal{A}_{t}$-measurable. Since addition, multiplication, and
exponentiation are continuous operations, every finite partial sum 
\begin{equation*}
F_{t,N}^{B}(x)=\sum_{\substack{ n\in B  \\ n\leq N}}a_{n}(t,x)\frac{%
g(t,x)^{n}}{n!},
\end{equation*}%
is $\mathcal{A}_{t}$-measurable. Taking the pointwise limit wherever the
series converges, we conclude that $F_{t}^{B}(x)$ is $\mathcal{A}_{t}$%
-measurable. Therefore, the process is adapted.

\noindent \textbf{Step 2. Integrability.}

By assumption, 
\begin{equation*}
\mathbb{E}[|F_{t}^{B}(x)|]<\infty ,
\end{equation*}%
for every $t\geq 0$. Thus, the process is integrable.

\noindent \textbf{Step 3. Conditional expectation property.}

By assumption, for every $0\leq s\leq t$, 
\begin{equation*}
\mathbb{E}[F_{t}^{B}(x)\mid \mathcal{A}_{s}]=F_{s}^{B}(x)\quad \text{a.s.}
\end{equation*}%
This is precisely the martingale conditional expectation property.

Combining Steps 1--3, we conclude that 
\begin{equation*}
(F_{t}^{B}(x))_{t\geq 0},
\end{equation*}%
is a martingale with respect to $\mathbb{F}$.
\end{proof}

The preceding theorem gives a direct characterization, but its third
assumption essentially states the martingale property itself. We therefore
provide a more useful sufficient condition based on the individual Taylor
terms.

Define 
\begin{equation*}
Y_{n,t}(x)=a_{n}(t,x)\frac{g(t,x)^{n}}{n!},
\end{equation*}%
such that%
\begin{equation*}
F_{t}^{B}(x)=\sum_{n\in B}Y_{n,t}(x).
\end{equation*}

\begin{theorem}[Martingale Taylor Terms]
Fix $x\in\mathbb{M}^p$ and $B\in\mathcal{B}(\mathbb{N}_0)$.

Suppose that:

\begin{enumerate}
\item For every $n\in B$, the process 
\begin{equation*}
(Y_{n,t}(x))_{t\geq 0},
\end{equation*}%
is an integrable martingale with respect to $\mathbb{F}$.

\item For every $t\geq 0$, 
\begin{equation*}
\mathbb{E}\left[ \sum_{n\in B}|Y_{n,t}(x)|\right] <\infty .
\end{equation*}%
Then 
\begin{equation*}
F_{t}^{B}(x)=\sum_{n\in B}Y_{n,t}(x),
\end{equation*}%
is a martingale with respect to $\mathbb{F}$.
\end{enumerate}
\end{theorem}

\begin{proof}
We verify the martingale properties.

\noindent \textbf{Step 1. Adaptedness.}

Since every $Y_{n,t}(x)$ is an adapted process, each finite partial sum 
\begin{equation*}
F_{t,N}^{B}(x)=\sum_{\substack{ n\in B  \\ n\leq N}}Y_{n,t}(x),
\end{equation*}%
is $\mathcal{A}_{t}$-measurable. By the assumed absolute convergence, 
\begin{equation*}
F_{t}^{B}(x)=\lim_{N\rightarrow \infty }F_{t,N}^{B}(x),
\end{equation*}%
almost surely. Hence, $F_{t}^{B}(x)$ is measurable with respect to the
completed sigma-algebra $\mathcal{A}_{t}$, or admits an $\mathcal{A}_{t}$%
-measurable version.

\noindent \textbf{Step 2. Integrability.}

Using the triangle inequality, 
\begin{equation*}
\left\vert F_{t}^{B}(x)\right\vert \leq \sum_{n\in B}|Y_{n,t}(x)|.
\end{equation*}%
Taking expectations gives 
\begin{equation*}
\mathbb{E}[|F_{t}^{B}(x)|]\leq \mathbb{E}\left[ \sum_{n\in B}|Y_{n,t}(x)|%
\right] <\infty .
\end{equation*}%
Therefore, $F_{t}^{B}(x)$ is integrable.

\noindent \textbf{Step 3. Conditional expectation.}

Fix $0\leq s\leq t$. Since every $Y_{n,t}(x)$ is a martingale, 
\begin{equation*}
\mathbb{E}[Y_{n,t}(x)\mid \mathcal{A}_{s}]=Y_{n,s}(x).
\end{equation*}%
By the assumed integrability, 
\begin{equation*}
\mathbb{E}\left[ \sum_{n\in B}|Y_{n,t}(x)|\right] <\infty .
\end{equation*}%
Consequently, the series 
\begin{equation*}
\sum_{n\in B}Y_{n,t}(x),
\end{equation*}%
converges in $L^{1}$. Conditional expectation is a continuous linear
operator on $L^{1}$. Therefore, we can write 
\begin{align*}
\mathbb{E}[F_{t}^{B}(x)\mid \mathcal{A}_{s}]& =\mathbb{E}\left[ \sum_{n\in
B}Y_{n,t}(x)\middle|\mathcal{A}_{s}\right] =\sum_{n\in B}\mathbb{E}%
[Y_{n,t}(x)\mid \mathcal{A}_{s}] \\
& =\sum_{n\in B}Y_{n,s}(x)=F_{s}^{B}(x),
\end{align*}%
and thus 
\begin{equation*}
\mathbb{E}[F_{t}^{B}(x)\mid \mathcal{A}_{s}]=F_{s}^{B}(x)\quad \text{a.s.}
\end{equation*}%
All three martingale conditions are satisfied and hence $F_{t}^{B}(x)$ is a
martingale.
\end{proof}

\begin{example}[A Simple Example]
Consider the special case in which the Taylor parameter is deterministic and
constant, i.e., 
\begin{equation*}
g(t,\omega ,x)=c,
\end{equation*}%
where $c\in \mathbb{M}$. Suppose that, for every $n\in \mathbb{N}_{0}$, the
process%
\begin{equation*}
(a_{n}(t))_{t\geq 0},
\end{equation*}%
is an integrable martingale. Then 
\begin{equation*}
Y_{n,t}=a_{n}(t)\frac{c^{n}}{n!},
\end{equation*}%
is also a martingale because multiplication by a deterministic constant
preserves the martingale property. If 
\begin{equation*}
\mathbb{E}\left[ \sum_{n=0}^{\infty }|a_{n}(t)|\frac{|c|^{n}}{n!}\right]
<\infty ,
\end{equation*}%
for every $t\geq 0$, then 
\begin{equation*}
F_{t}=\sum_{n=0}^{\infty }a_{n}(t)\frac{c^{n}}{n!},
\end{equation*}%
is a martingale. This provides a direct example of a SMTMF that is a
martingale.
\end{example}

It is important to distinguish between the following statements:

\begin{enumerate}
\item Every Taylor coefficient process $a_n(t)$ is a martingale.

\item Every Taylor term 
\begin{equation*}
a_{n}(t)\frac{g(t)^{n}}{n!},
\end{equation*}%
is a martingale.

\item The complete Taylor series 
\begin{equation*}
F_{t}^{B}=\sum_{n\in B}a_{n}(t)\frac{g(t)^{n}}{n!},
\end{equation*}%
is a martingale.
\end{enumerate}

The first condition does not automatically imply the second if $g(t)$ is
random. For example, the product of two martingales is not necessarily a
martingale. Therefore, when the Taylor parameter is random, it is important
to impose conditions directly on the Taylor terms or establish their
conditional expectation properties.

\begin{remark}
The martingale property is a temporal property. Consequently, a stochastic
multivariate Taylor measure function is a martingale only after a time
parameter and an appropriate filtration have been specified. For a random
function indexed by $x$, the martingale property is generally understood
pointwise in $x$. If one additionally requires the martingale property
simultaneously for all spatial inputs, further regularity and separability
assumptions may be needed.
\end{remark}

A SMTMF is a martingale when its associated process satisfies adaptedness,
integrability, and the conditional expectation property. A particularly
useful sufficient condition is that each individual Taylor term 
\begin{equation*}
a_{n}(t,x)\frac{g(t,x)^{n}}{n!},
\end{equation*}%
is a martingale and that the series is sufficiently integrable to permit
interchanging summation and conditional expectation. This provides a natural
foundation for investigating martingale valued stochastic Taylor measure
functions and their applications in stochastic processes and random fields.

\subsection{Martingale Conditions on the SMTMF parameter and coefficients}

Let $(\Omega ,\mathcal{A},\mathbb{F},P)$ be a filtered probability space,
where 
\begin{equation*}
\mathbb{F}=(\mathcal{A}_{t})_{t\geq 0},
\end{equation*}%
is a filtration satisfying the usual conditions.

Let $p\in \mathbb{N}^{+}$ and $\mathbb{M}\in \{\mathbb{R},\mathbb{C}\}$.
Consider a stochastic multivariate Taylor measure function of the form 
\begin{equation*}
F_{t}^{B}(x)=\sum_{n\in B}a_{n}(t,\omega ,x)\frac{g(t,\omega ,x)^{n}}{n!},
\end{equation*}%
where $B\in \mathcal{B}(\mathbb{N}_{0})$. For every fixed $x\in \mathbb{M}%
^{p}$, write 
\begin{equation*}
A_{n,t}=a_{n}(t,\cdot ,x),\qquad G_{t}=g(t,\cdot ,x),
\end{equation*}%
such that 
\begin{equation*}
F_{t}^{B}(x)=\sum_{n\in B}Y_{n,t},
\end{equation*}%
where 
\begin{equation*}
Y_{n,t}=\frac{A_{n,t}G_{t}^{n}}{n!}.
\end{equation*}%
Our objective is to establish sufficient conditions on the coefficient
processes $A_{n,t}$ and the random Taylor parameter $G_{t}$ that guarantee
that every individual Taylor term is a martingale.

We begin with conditional covariance conditions.

\begin{definition}[Conditional Covariance]
For integrable random variables $U,V$ such that $UV$ is integrable, define 
\begin{equation*}
\func{Cov}(U,V\mid \mathcal{A}_{s})=\mathbb{E}[UV\mid \mathcal{A}_{s}]-%
\mathbb{E}[U\mid \mathcal{A}_{s}]\mathbb{E}[V\mid \mathcal{A}_{s}].
\end{equation*}%
For complex-valued variables, we use the bilinear convention without complex
conjugation.
\end{definition}

The criterion is proved next.

\begin{theorem}[Conditional Covariance Criterion]
Fix $n\in \mathbb{N}_{0}$. Suppose that the processes $(A_{n,t})_{t\geq 0}$
and $(G_{t}^{n})_{t\geq 0}$ are integrable martingales with respect to $%
\mathbb{F}$. Assume also that 
\begin{equation*}
\mathbb{E}[|A_{n,t}G_{t}^{n}|]<\infty ,
\end{equation*}%
for every $t\geq 0$. Then 
\begin{equation*}
Y_{n,t}=\frac{A_{n,t}G_{t}^{n}}{n!},
\end{equation*}%
is a martingale if and only if 
\begin{equation*}
{\func{Cov}(A_{n,t},G_{t}^{n}\mid \mathcal{A}_{s})=0},
\end{equation*}%
for every $0\leq s\leq t$.
\end{theorem}

\begin{proof}
Fix $0\leq s\leq t$. By the definition of conditional covariance, 
\begin{equation*}
\mathbb{E}[A_{n,t}G_{t}^{n}\mid \mathcal{A}_{s}]=\mathbb{E}[A_{n,t}\mid 
\mathcal{A}_{s}]\mathbb{E}[G_{t}^{n}\mid \mathcal{A}_{s}]+\func{Cov}%
(A_{n,t},G_{t}^{n}\mid \mathcal{A}_{s}).
\end{equation*}%
Since $A_{n,t}$ is a martingale, 
\begin{equation*}
\mathbb{E}[A_{n,t}\mid \mathcal{A}_{s}]=A_{n,s}.
\end{equation*}%
Similarly, since $G_{t}^{n}$ is a martingale, 
\begin{equation*}
\mathbb{E}[G_{t}^{n}\mid \mathcal{A}_{s}]=G_{s}^{n}.
\end{equation*}%
Consequently, 
\begin{equation*}
\mathbb{E}[Y_{n,t}\mid \mathcal{A}_{s}]=\frac{A_{n,s}G_{s}^{n}}{n!}+\frac{1}{%
n!}\func{Cov}(A_{n,t},G_{t}^{n}\mid \mathcal{A}_{s}).
\end{equation*}%
Therefore, 
\begin{equation*}
\mathbb{E}[Y_{n,t}\mid \mathcal{A}_{s}]=Y_{n,s}
\end{equation*}%
if and only if 
\begin{equation*}
\func{Cov}(A_{n,t},G_{t}^{n}\mid \mathcal{A}_{s})=0.
\end{equation*}%
Since the product is adapted and integrable, this proves the result.
\end{proof}

Conditional independence of $A_{n,t}$ and $G_{t}^{n}$ given $\mathcal{A}_{s}$
is sufficient for the conditional covariance to vanish, provided the
relevant integrability conditions hold. However, we note that conditional
independence is stronger than necessary. Vanishing conditional covariance is
sufficient.

\subsection{Semimartingale Characterization}

We now obtain an explicit condition involving the drift and volatility of
the coefficient processes and the random Taylor parameter. Throughout this
section, assume that all semimartingales are real-valued and continuous.
Suppose that 
\begin{equation*}
dA_{n,t}=\alpha _{n,t}dt+dM_{n,t}
\end{equation*}%
and 
\begin{equation*}
dG_{t}=\beta _{t}dt+dN_{t},
\end{equation*}%
where $M_{n,t}$ and $N_{t}$ are continuous local martingales. Assume that 
\begin{equation*}
d\langle M_{n},N\rangle _{t}=q_{n,t}\,dt,
\end{equation*}%
and 
\begin{equation*}
d\langle N\rangle _{t}=v_{t}\,dt.
\end{equation*}%
Here, $\alpha _{n,t}$, $\beta _{t}$, $q_{n,t}$, and $v_{t}$ are predictable
processes.

\begin{theorem}[Semimartingale Martingale Criterion]
\label{semimgthm}Let $n\in \mathbb{N}^{+}$ and define 
\begin{equation*}
Y_{n,t}=\frac{A_{n,t}G_{t}^{n}}{n!}.
\end{equation*}%
Suppose that $A_{n,t}$ and $G_{t}$ satisfy the preceding semimartingale
decompositions. If 
\begin{equation}
0=G_{t}^{n}\alpha _{n,t}+nA_{n,t}G_{t}^{n-1}\beta _{t}+nG_{t}^{n-1}q_{n,t}+%
\frac{n(n-1)}{2}A_{n,t}G_{t}^{n-2}v_{t},  \label{MGcriterion}
\end{equation}%
for $dt\otimes \mathbb{P}$-almost every $(t,\omega )$, then $Y_{n,t}$ is a
local martingale.

If, in addition, $Y_{n,t}$ is of class $D$ on every finite time interval,
then it is a true martingale. For $n=0$, the corresponding condition is 
\begin{equation*}
\alpha _{0,t}=0.
\end{equation*}
\end{theorem}

\begin{proof}
We use It\^{o}'s formula (see \cite{KaratzasShreve1991}, \cite{RevuzYor1999}%
, \cite{RogersWilliams2000}, \cite{micheas2018theory}) and the
semimartingale product rule.

\noindent \textbf{Step 1. Apply It\^{o}'s formula to }$G_{t}^{n}$

Consider 
\begin{equation*}
f(z)=z^{n}.
\end{equation*}%
We have trivially 
\begin{equation*}
f^{\prime }(z)=nz^{n-1}
\end{equation*}%
and 
\begin{equation*}
f^{\prime \prime }(z)=n(n-1)z^{n-2}.
\end{equation*}%
By It\^{o}'s formula, 
\begin{equation*}
d(G_{t}^{n})=nG_{t}^{n-1}\,dG_{t}+\frac{n(n-1)}{2}G_{t}^{n-2}\,d\langle
N\rangle _{t}.
\end{equation*}%
Substituting 
\begin{equation*}
dG_{t}=\beta _{t}\,dt+dN_{t}
\end{equation*}%
and 
\begin{equation*}
d\langle N\rangle _{t}=v_{t}\,dt,
\end{equation*}%
we obtain 
\begin{equation*}
d(G_{t}^{n})=nG_{t}^{n-1}\beta _{t}\,dt+nG_{t}^{n-1}\,dN_{t}+\frac{n(n-1)}{2}%
G_{t}^{n-2}v_{t}\,dt.
\end{equation*}

\noindent \textbf{Step 2. Apply the product rule}

Recall that 
\begin{equation*}
Y_{n,t}=\frac{1}{n!}A_{n,t}G_{t}^{n}.
\end{equation*}%
By integration by parts, 
\begin{equation*}
d(A_{n,t}G_{t}^{n})=G_{t}^{n}\,dA_{n,t}+A_{n,t}\,d(G_{t}^{n})+d\langle
M_{n},G^{n}\rangle _{t}.
\end{equation*}%
The quadratic covariation is 
\begin{equation*}
d\langle M_{n},G^{n}\rangle _{t}=nG_{t}^{n-1}d\langle M_{n},N\rangle _{t},
\end{equation*}%
and hence 
\begin{equation*}
d\langle M_{n},G^{n}\rangle _{t}=nG_{t}^{n-1}q_{n,t}\,dt.
\end{equation*}

\noindent \textbf{Step 3. Substitute the decompositions}

Using 
\begin{equation*}
dA_{n,t}=\alpha _{n,t}\,dt+dM_{n,t},
\end{equation*}%
we obtain 
\begin{align*}
d(A_{n,t}G_{t}^{n})& =G_{t}^{n}(\alpha _{n,t}\,dt+dM_{n,t})+A_{n,t}\Bigg[%
nG_{t}^{n-1}\beta _{t}\,dt+nG_{t}^{n-1}\,dN_{t} \\
& +\frac{n(n-1)}{2}G_{t}^{n-2}v_{t}\,dt\Bigg]+nG_{t}^{n-1}q_{n,t}\,dt.
\end{align*}%
Collecting the drift terms gives 
\begin{equation*}
d(A_{n,t}G_{t}^{n})=D_{n,t}\,dt+G_{t}^{n}\,dM_{n,t}+nA_{n,t}G_{t}^{n-1}%
\,dN_{t},
\end{equation*}%
where 
\begin{equation*}
D_{n,t}=G_{t}^{n}\alpha _{n,t}+nA_{n,t}G_{t}^{n-1}\beta
_{t}+nG_{t}^{n-1}q_{n,t}+\frac{n(n-1)}{2}A_{n,t}G_{t}^{n-2}v_{t}.
\end{equation*}%
Dividing by $n!$ yields 
\begin{equation*}
dY_{n,t}=\frac{1}{n!}\left[ D_{n,t}\,dt+G_{t}^{n}%
\,dM_{n,t}+nA_{n,t}G_{t}^{n-1}\,dN_{t}\right] .
\end{equation*}

\noindent \textbf{Step 4. Eliminate the drift}

By hypothesis, we have 
\begin{equation*}
D_{n,t}=0,
\end{equation*}%
for $dt\otimes \mathbb{P}$-almost every $(t,\omega )$. Therefore, we can
write%
\begin{equation*}
dY_{n,t}=\frac{1}{n!}\left[ G_{t}^{n}\,dM_{n,t}+nA_{n,t}G_{t}^{n-1}\,dN_{t}%
\right] .
\end{equation*}%
The right-hand side is a stochastic integral with respect to continuous
local martingales. Consequently, $Y_{n,t}$ is a local martingale.

\noindent \textbf{Step 5. Establish the true martingale property}

A local martingale need not be a true martingale. However, a local
martingale of class $D$ on a finite time interval is a uniformly integrable
martingale on that interval. Thus, if $Y_{n,t}$ is of class $D$ on every
finite time interval, then 
\begin{equation*}
\mathbb{E}[Y_{n,t}\mid \mathcal{A}_{s}]=Y_{n,s},
\end{equation*}%
for every $0\leq s\leq t$. This completes the proof.
\end{proof}

\subsection{Brownian Motion Specialization}

Suppose that $W_{t}$ is a $d$-dimensional Brownian motion and 
\begin{equation*}
dA_{n,t}=\alpha _{n,t}\,dt+\boldsymbol{\sigma}_{n,t}\cdot dW_{t},
\end{equation*}%
and 
\begin{equation*}
dG_{t}=\beta _{t}\,dt+\boldsymbol{\tau}_{t}\cdot dW_{t}.
\end{equation*}

Here, 
\begin{equation*}
\boldsymbol{\sigma}_{n,t}=(\sigma _{n,t}^{1},\ldots ,\sigma _{n,t}^{d}),
\end{equation*}%
and 
\begin{equation*}
\boldsymbol{\tau}_{t}=(\tau _{t}^{1},\ldots ,\tau _{t}^{d}).
\end{equation*}

The quadratic covariations are 
\begin{equation*}
d\langle M_{n},N\rangle _{t}=(\boldsymbol{\sigma}_{n,t}\cdot %
\boldsymbol{\tau}_{t})\,dt,
\end{equation*}%
and 
\begin{equation*}
d\langle N\rangle _{t}=\Vert \boldsymbol{\tau}_{t}\Vert ^{2}\,dt.
\end{equation*}

\begin{corollary}[Brownian Martingale Condition]
Under the preceding assumptions, 
\begin{equation*}
Y_{n,t}=\frac{A_{n,t}G_{t}^{n}}{n!},
\end{equation*}%
is a local martingale if 
\begin{equation*}
0=G_{t}^{n}\alpha _{n,t}+nA_{n,t}G_{t}^{n-1}\beta _{t}+nG_{t}^{n-1}(%
\boldsymbol{\sigma}_{n,t}\cdot \boldsymbol{\tau}_{t})+\frac{n(n-1)}{2}%
A_{n,t}G_{t}^{n-2}\Vert \boldsymbol{\tau}_{t}\Vert ^{2}.
\end{equation*}

If it is also of class $D$ on every finite time interval, it is a true
martingale.
\end{corollary}

\begin{proof}
Substitute 
\begin{equation*}
q_{n,t}=\boldsymbol{\sigma}_{n,t}\cdot \boldsymbol{\tau}_{t}
\end{equation*}%
and 
\begin{equation*}
v_{t}=\Vert \boldsymbol{\tau}_{t}\Vert ^{2}
\end{equation*}%
into Theorem \ref{semimgthm}. The conclusion follows immediately.
\end{proof}

\subsection{Explicit Examples}

\begin{example}[Deterministic Taylor Parameter]
Suppose that 
\begin{equation*}
G_{t}=c,
\end{equation*}%
where $c\neq 0$ is deterministic. Then 
\begin{equation*}
\beta _{t}=0,\qquad \boldsymbol{\tau}_{t}=0.
\end{equation*}%
The Brownian martingale condition reduces to 
\begin{equation*}
c^{n}\alpha _{n,t}=0,
\end{equation*}%
and therefore 
\begin{equation*}
\alpha _{n,t}=0.
\end{equation*}%
In particular, if each $A_{n,t}$ is a true martingale and 
\begin{equation*}
\mathbb{E}\left[ \sum_{n\in B}|A_{n,t}|\frac{|c|^{n}}{n!}\right] <\infty ,
\end{equation*}%
then 
\begin{equation*}
F_{t}^{B}=\sum_{n\in B}A_{n,t}\frac{c^{n}}{n!},
\end{equation*}%
is a true martingale.
\end{example}

\begin{example}[Correlated Brownian-Driven Processes]
Consider 
\begin{equation*}
dA_{n,t}=\alpha _{n,t}\,dt+\sigma _{n,t}\,dW_{t},
\end{equation*}%
and 
\begin{equation*}
dG_{t}=\beta _{t}\,dt+\tau _{t}\,dW_{t}.
\end{equation*}%
The martingale condition is 
\begin{equation*}
0=G_{t}^{n}\alpha _{n,t}+nA_{n,t}G_{t}^{n-1}\beta _{t}+nG_{t}^{n-1}\sigma
_{n,t}\tau _{t}+\frac{n(n-1)}{2}A_{n,t}G_{t}^{n-2}\tau _{t}^{2}.
\end{equation*}%
For $n=1$, this becomes 
\begin{equation*}
G_{t}\alpha _{1,t}+A_{1,t}\beta _{t}+\sigma _{1,t}\tau _{t}=0.
\end{equation*}%
For $n=2$, we obtain 
\begin{equation*}
0=G_{t}^{2}\alpha _{2,t}+2A_{2,t}G_{t}\beta _{t}+2G_{t}\sigma _{2,t}\tau
_{t}+A_{2,t}\tau _{t}^{2}.
\end{equation*}%
These equations explicitly demonstrate how the drift and volatility of the
coefficient processes must compensate for the drift, volatility, and
quadratic variation of the random Taylor parameter.
\end{example}

\begin{example}[Exponential Brownian Taylor Parameter]
Let 
\begin{equation*}
G_{t}=\exp \left( W_{t}-\frac{t}{2}\right) ,
\end{equation*}%
where $W_{t}$ is standard Brownian motion. It\^{o}'s formula gives 
\begin{equation*}
dG_{t}=G_{t}dW_{t}.
\end{equation*}%
Therefore, we have%
\begin{equation*}
\beta _{t}=0,\qquad \tau _{t}=G_{t}.
\end{equation*}%
For every $n\geq 1$, 
\begin{equation*}
d(G_{t}^{n})=nG_{t}^{n}\,dW_{t}+\frac{n(n-1)}{2}G_{t}^{n}dt.
\end{equation*}%
Thus, $G_{t}$ is a martingale, but $G_{t}^{n}$ is not a martingale for $%
n\geq 2$. Now choose 
\begin{equation*}
dA_{n,t}=-\frac{n(n-1)}{2}A_{n,t}dt+dM_{n,t},
\end{equation*}%
where $M_{n,t}$ is a continuous local martingale orthogonal to $W_{t}$. Then 
\begin{equation*}
q_{n,t}=0,
\end{equation*}%
and 
\begin{equation*}
\alpha _{n,t}=-\frac{n(n-1)}{2}A_{n,t}.
\end{equation*}%
The drift condition becomes 
\begin{equation*}
D_{n,t}=-\frac{n(n-1)}{2}A_{n,t}G_{t}^{n}+\frac{n(n-1)}{2}A_{n,t}G_{t}^{n}=0.
\end{equation*}%
Consequently%
\begin{equation*}
Y_{n,t}=\frac{A_{n,t}G_{t}^{n}}{n!},
\end{equation*}%
is a local martingale. If it is of class $D$, it is a true martingale.
\end{example}

\subsection{Martingale Property of the Complete STMTF}

\begin{theorem}[Martingale Property of the Complete SMTMF]
Let $B\in \mathcal{B}(\mathbb{N}_{0})$. Suppose that, for every $n\in B$, 
\begin{equation*}
Y_{n,t}=\frac{A_{n,t}G_{t}^{n}}{n!},
\end{equation*}%
is a true martingale with respect to $\mathbb{F}$. Assume also that 
\begin{equation*}
{\mathbb{E}\left[ \sum_{n\in B}\left\vert A_{n,t}\frac{G_{t}^{n}}{n!}%
\right\vert \right] <\infty ,}
\end{equation*}%
for every $t\geq 0$. Then 
\begin{equation*}
F_{t}^{B}=\sum_{n\in B}A_{n,t}\frac{G_{t}^{n}}{n!},
\end{equation*}%
is a true martingale with respect to $\mathbb{F}$.
\end{theorem}

\begin{proof}
First, since every $Y_{n,t}$ is adapted, each finite partial sum 
\begin{equation*}
F_{t,N}^{B}=\sum_{\substack{ n\in B  \\ n\leq N}}Y_{n,t},
\end{equation*}%
is $\mathcal{A}_{t}$-measurable. The assumed absolute convergence ensures
that the limit 
\begin{equation*}
F_{t}^{B}=\lim_{N\rightarrow \infty }F_{t,N}^{B},
\end{equation*}%
exists almost surely and is measurable with respect to the completed
filtration. Second, by the triangle inequality, 
\begin{equation*}
\left\vert F_{t}^{B}\right\vert \leq \sum_{n\in B}|Y_{n,t}|.
\end{equation*}%
Hence, we have 
\begin{equation*}
\mathbb{E}[|F_{t}^{B}|]\leq \mathbb{E}\left[ \sum_{n\in B}|Y_{n,t}|\right]
<\infty .
\end{equation*}%
Thus, $F_{t}^{B}$ is integrable.\newline
Finally, fix $0\leq s\leq t$. The assumed integrability implies convergence
of the partial sums in $L^{1}$. Since conditional expectation is continuous
in $L^{1}$, 
\begin{equation*}
\mathbb{E}[F_{t}^{B}\mid \mathcal{A}_{s}]=\mathbb{E}\left[ \sum_{n\in
B}Y_{n,t}\middle|\mathcal{A}_{s}\right] =\sum_{n\in B}\mathbb{E}[Y_{n,t}\mid 
\mathcal{A}_{s}]=\sum_{n\in B}Y_{n,s}=F_{s}^{B}.
\end{equation*}%
Therefore, $F_{t}^{B}$ is a true martingale.
\end{proof}

We have established conditions on the coefficient processes and random
Taylor parameter that guarantee the martingale property of individual Taylor
terms. The central semimartingale condition of equation (\ref{MGcriterion})
reveals four contributions:

\begin{enumerate}
\item The drift of the random Taylor coefficient.

\item The drift of the random Taylor parameter.

\item The quadratic covariation between the coefficient and the Taylor
parameter.

\item The quadratic variation of the Taylor parameter.
\end{enumerate}

In particular, the coefficient processes do not need to be martingales
individually. Their drift components may compensate for the drift introduced
by the random Taylor parameter. The resulting theory provides a
semimartingale framework for studying martingale valued stochastic
multivariate Taylor measure functions.

\section{Explicit Construction of Brownian motion via SMTMFs}

In order to exemplify the flexibility of SMTMFs, we present in explicit
details how we can construct some of the most important stochastic
processes, and we begin with Brownian motion.

The construction is based on the L\'{e}vy-Ciesielski representation of
Brownian motion (see \cite{Ciesielski1961}, \cite{Kuo1975}, \cite%
{Lifshits2012}), which expresses Brownian motion as an almost surely
uniformly convergent series of independent standard Gaussian random
variables multiplied by deterministic Schauder functions.

We show that every finite approximation in this representation is a SMTMF
associated with a finite Taylor measure. We then prove that these finite
approximations converge almost surely uniformly to a process with continuous
sample paths, centered Gaussian finite dimensional distributions, and
covariance function $\min (s,t)$. Consequently, the limiting process is
standard Brownian motion.

Let $(\Omega ,\mathcal{A},\mathbb{P})$ be a probability space. Recall that a
stochastic MTMF has the form 
\begin{equation*}
F^{B}(\omega ,x)=\sum_{n\in B}a_{n}(\omega ,x)\frac{g(\omega ,x)^{n}}{n!},
\end{equation*}%
where $B\in \mathcal{B}(\mathbb{N}_{0})$ and the coefficient functions and
Taylor parameter satisfy the required measurability and finiteness
conditions. For our construction, we take 
\begin{equation*}
p=1,\qquad x=t\in \lbrack 0,1],
\end{equation*}%
and choose the deterministic Taylor parameter 
\begin{equation*}
g(t)=1.
\end{equation*}%
We construct finite stochastic MTMFs whose almost sure uniform limit is
standard Brownian motion.

\begin{definition}[Schauder Functions]
For $j\in \mathbb{N}_{0}$ and $k\in \{0,1,\ldots ,2^{j}-1\}$, define 
\begin{equation*}
h_{j,k}(t)=%
\begin{cases}
2^{j+1}t-2k, & \frac{k}{2^{j}}\leq t\leq \frac{2k+1}{2^{j+1}}, \\[6pt] 
2k+2-2^{j+1}t, & \frac{2k+1}{2^{j+1}}\leq t\leq \frac{k+1}{2^{j}}, \\[6pt] 
0, & \text{otherwise}.%
\end{cases}%
\end{equation*}%
These functions are continuous triangular functions supported on 
\begin{equation*}
\left[ \frac{k}{2^{j}},\frac{k+1}{2^{j}}\right] ,
\end{equation*}%
and they satisfy 
\begin{equation*}
0\leq h_{j,k}(t)\leq 1.
\end{equation*}
\end{definition}

Let 
\begin{equation*}
Z_{\ast },\quad Z_{j,k},
\end{equation*}%
where $j\in \mathbb{N}_{0}$ and $0\leq k<2^{j}$, be mutually independent
standard normal random variables. i.e., 
\begin{equation*}
Z_{\ast },Z_{j,k}\sim N(0,1).
\end{equation*}%
Define 
\begin{equation*}
c_{j}=2^{-(j+2)/2}.
\end{equation*}%
For each $N\in \mathbb{N}_{0}$, define the finite index set 
\begin{equation*}
B_{N}=\{0\}\cup \left\{ 2^{j}+k:0\leq j\leq N,\;0\leq k<2^{j}\right\} .
\end{equation*}%
Define the coefficient functions by 
\begin{equation*}
a_{0}(\omega ,t)=Z_{\ast }(\omega )t,
\end{equation*}%
and 
\begin{equation*}
a_{2^{j}+k}(\omega ,t)=(2^{j}+k)!\,Z_{j,k}(\omega )c_{j}h_{j,k}(t).
\end{equation*}%
Finally, let 
\begin{equation*}
g(\omega ,t)=1.
\end{equation*}%
The resulting stochastic MTMF is 
\begin{equation*}
F_{N}(\omega ,t)=\sum_{n\in B_{N}}a_{n}(\omega ,t)\frac{g(\omega ,t)^{n}}{n!}%
,
\end{equation*}%
or equivalently, 
\begin{equation}
F_{N}(t)=Z_{\ast
}t+\sum_{j=0}^{N}\sum_{k=0}^{2^{j}-1}Z_{j,k}2^{-(j+2)/2}h_{j,k}(t).  \tag{1}
\label{F_nt}
\end{equation}

\begin{theorem}[Brownian Motion Construction Using Stochastic MTMFs]
Let $(\Omega ,\mathcal{A},\mathbb{P})$ support mutually independent standard
normal random variables 
\begin{equation*}
Z_{\ast },\quad Z_{j,k},\qquad j\in \mathbb{N}_{0},\quad 0\leq k<2^{j}.
\end{equation*}%
For every $N\in \mathbb{N}_{0}$, define $F_{N}$ by equation (\ref{F_nt}).

Then the following statements hold.

\begin{enumerate}
\item Every $F_N$ is a finite stochastic multivariate Taylor measure
function.

\item The sequence $(F_N)$ converges almost surely uniformly on $[0,1]$.

\item The limiting process 
\begin{equation*}
W_{t}=\lim_{N\rightarrow \infty }F_{N}(t),
\end{equation*}%
has continuous sample paths and satisfies 
\begin{equation*}
W_{0}=0.
\end{equation*}

\item The limiting process is centered Gaussian and has covariance 
\begin{equation*}
{\mathbb{E}[W_{s}W_{t}]=\min (s,t).}
\end{equation*}

\item The process $W$ is standard Brownian motion with respect to its
natural filtration.
\end{enumerate}
\end{theorem}

\begin{proof}
\textbf{Proof of statement 1: Finite SMTMF Representation}

By construction, 
\begin{equation*}
B_{N}=\{0\}\cup \left\{ 2^{j}+k:0\leq j\leq N,\;0\leq k<2^{j}\right\} ,
\end{equation*}%
is finite and%
\begin{equation*}
g(t)=1.
\end{equation*}%
The coefficient functions are 
\begin{equation*}
a_{0}(\omega ,t)=Z_{\ast }(\omega )t
\end{equation*}%
and 
\begin{equation*}
a_{2^{j}+k}(\omega ,t)=(2^{j}+k)!Z_{j,k}(\omega )2^{-(j+2)/2}h_{j,k}(t).
\end{equation*}%
Each coefficient is measurable in $(\omega ,t)$ because the Gaussian random
variables are measurable and the Schauder functions are continuous.
Therefore, we have that%
\begin{equation*}
F_{N}(\omega ,t)=\sum_{n\in B_{N}}a_{n}(\omega ,t)\frac{1}{n!},
\end{equation*}%
is a finite stochastic MTMF.

\textbf{Proof of statement 2: Uniform Convergence}

For every $j\in \mathbb{N}_{0}$, define 
\begin{equation*}
M_{j}=\max_{0\leq k<2^{j}}|Z_{j,k}|.
\end{equation*}%
Since the Schauder functions at a fixed level have disjoint interiors and
satisfy 
\begin{equation*}
0\leq h_{j,k}(t)\leq 1,
\end{equation*}%
we have 
\begin{equation*}
\sup_{t\in \lbrack 0,1]}\left\vert
\sum_{k=0}^{2^{j}-1}Z_{j,k}c_{j}h_{j,k}(t)\right\vert \leq c_{j}M_{j}.
\end{equation*}%
We will prove that 
\begin{equation*}
\sum_{j=0}^{\infty }c_{j}M_{j}<\infty \quad \text{almost surely}.
\end{equation*}%
For a standard normal random variable $Z$, the Gaussian tail estimate gives 
\begin{equation*}
\mathbb{P}(|Z|>u)\leq 2e^{-u^{2}/2}.
\end{equation*}%
Using the union bound, 
\begin{equation*}
\mathbb{P}(M_{j}>u)\leq \sum_{k=0}^{2^{j}-1}\mathbb{P}(|Z_{j,k}|>u)\leq
2^{j+1}e^{-u^{2}/2}.
\end{equation*}%
Choose $C>\sqrt{2\log 2}$ and set 
\begin{equation*}
u_{j}=C\sqrt{j+1}.
\end{equation*}%
Then 
\begin{equation*}
\mathbb{P}(M_{j}>u_{j})\leq 2^{j+1}e^{-C^{2}(j+1)/2}=2\exp \left[ -\left( 
\frac{C^{2}}{2}-\log 2\right) (j+1)\right] .
\end{equation*}%
Therefore, 
\begin{equation*}
\sum_{j=0}^{\infty }\mathbb{P}(M_{j}>u_{j})<\infty .
\end{equation*}%
By the first Borel--Cantelli lemma, almost surely there exists a finite
random integer $J$ such that 
\begin{equation*}
M_{j}\leq C\sqrt{j+1},
\end{equation*}%
for all $j\geq J$. Consequently, 
\begin{equation*}
\sum_{j=0}^{\infty }c_{j}M_{j}=\sum_{j=0}^{\infty }2^{-(j+2)/2}M_{j}<\infty .
\end{equation*}%
Indeed, the tail is bounded by 
\begin{equation*}
C\sum_{j=J}^{\infty }2^{-(j+2)/2}\sqrt{j+1},
\end{equation*}%
which converges. Now let $N>M$. Then we have 
\begin{equation*}
\Vert F_{N}-F_{M}\Vert _{\infty }=\sup_{t\in \lbrack 0,1]}\left\vert
F_{N}(t)-F_{M}(t)\right\vert \leq \sum_{j=M+1}^{N}c_{j}M_{j}.
\end{equation*}%
The right-hand side converges to zero almost surely as $M,N\rightarrow
\infty $. Hence, $(F_{N})$ is almost surely a Cauchy sequence in $C([0,1])$.
Since $C([0,1])$ is complete, there exists a continuous random function $W$
such that 
\begin{equation*}
\Vert F_{N}-W\Vert _{\infty }\longrightarrow 0\quad \text{almost surely},
\end{equation*}%
which proves the claim.

\textbf{Proof of statement 3: Continuity and Initial Value}

Every $F_{N}$ is a finite linear combination of continuous Schauder
functions and the continuous function $t\mapsto Z_{\ast }t$. Therefore, $%
F_{N}$ has continuous sample paths. By statement 2, $F_{N}$ converges
uniformly almost surely to $W$. The uniform limit of continuous functions is
continuous, and thus, $W$ has continuous sample paths almost surely.

Furthermore, 
\begin{equation*}
F_{N}(0)=0,
\end{equation*}%
for every $N$, because 
\begin{equation*}
h_{j,k}(0)=0.
\end{equation*}%
Taking limits, 
\begin{equation*}
{W_{0}=0,}
\end{equation*}%
which proves statement 3.

\textbf{Proof of statement 4: Gaussianity and Covariance}

By construction, 
\begin{equation*}
W_{t}=Z_{\ast }t+\sum_{j=0}^{\infty
}\sum_{k=0}^{2^{j}-1}Z_{j,k}c_{j}h_{j,k}(t).
\end{equation*}%
Since all Gaussian variables are centered and mutually independent, 
\begin{equation*}
\mathbb{E}[W_{t}]=0.
\end{equation*}%
For $s,t\in \lbrack 0,1]$, independence yields 
\begin{equation}
\mathbb{E}[W_{s}W_{t}]=st+\sum_{j=0}^{\infty
}\sum_{k=0}^{2^{j}-1}c_{j}^{2}h_{j,k}(s)h_{j,k}(t).  \label{covbrownianmot}
\end{equation}%
To evaluate this series, define the Haar functions 
\begin{equation*}
e_{\ast }(u)=1,
\end{equation*}%
and 
\begin{equation*}
e_{j,k}(u)=2^{j/2}\left( _{I_{j,k}^{L}}(u)-_{I_{j,k}^{R}}(u)\right) ,
\end{equation*}%
where 
\begin{equation*}
I_{j,k}^{L}=\left[ \frac{k}{2^{j}},\frac{2k+1}{2^{j+1}}\right) ,
\end{equation*}%
and 
\begin{equation*}
I_{j,k}^{R}=\left[ \frac{2k+1}{2^{j+1}},\frac{k+1}{2^{j}}\right) .
\end{equation*}%
The Haar functions form a complete orthonormal basis of $L^{2}([0,1])$.
Direct integration gives 
\begin{equation*}
\int_{0}^{t}e_{j,k}(u)du=2^{-(j+2)/2}h_{j,k}(t),
\end{equation*}%
and 
\begin{equation*}
\int_{0}^{t}e_{\ast }(u)du=t.
\end{equation*}%
Thus, equation (\ref{covbrownianmot}) becomes 
\begin{equation*}
\mathbb{E}[W_{s}W_{t}]=\left( \int_{0}^{s}e_{\ast }(u)\,du\right) \left(
\int_{0}^{t}e_{\ast }(u)\,du\right) +\sum_{j=0}^{\infty
}\sum_{k=0}^{2^{j}-1}\left( \int_{0}^{s}e_{j,k}(u)\,du\right) \left(
\int_{0}^{t}e_{j,k}(u)\,du\right) .
\end{equation*}%
By Parseval's identity, we have 
\begin{equation*}
\mathbb{E}[W_{s}W_{t}]=\left\langle \mathbf{1}_{[0,s]},\mathbf{1}%
_{[0,t]}\right\rangle _{L^{2}([0,1])}=\int_{0}^{1}\mathbf{1}_{[0,s]}(u)%
\mathbf{1}_{[0,t]}(u)\,du=\int_{0}^{\min (s,t)}1du=\min (s,t).
\end{equation*}%
Therefore, it follows that 
\begin{equation*}
E[W_{s}W_{t}]=\min (s,t).
\end{equation*}%
We next prove Gaussianity.

For arbitrary $t_{1},\ldots ,t_{m}\in \lbrack 0,1]$ and $\lambda _{1},\ldots
,\lambda _{m}\in \mathbb{R}$, consider 
\begin{equation*}
X_{N}=\sum_{r=1}^{m}\lambda _{r}F_{N}(t_{r}).
\end{equation*}%
Every $X_{N}$ is normally distributed, because it is a finite linear
combination of independent Gaussian random variables. Moreover, 
\begin{equation*}
X_{N}\longrightarrow X=\sum_{r=1}^{m}\lambda _{r}W_{t_{r}},
\end{equation*}%
almost surely. The variances converge to 
\begin{equation*}
\sigma ^{2}=\sum_{r=1}^{m}\sum_{\ell =1}^{m}\lambda _{r}\lambda _{\ell }\min
(t_{r},t_{\ell }).
\end{equation*}%
Consequently, the characteristic functions converge to 
\begin{equation*}
\mathbb{E}[e^{iuX}]=e^{-u^{2}\sigma ^{2}/2},
\end{equation*}%
and thus, $X$ is normally distributed. Since this holds for all choices of $%
\lambda _{1},\ldots ,\lambda _{m}$, the process $W$ is Gaussian, which
proves statement 4.

\textbf{Proof of statement 5: Independent Increments}

Let 
\begin{equation*}
0\leq s<t\leq 1.
\end{equation*}%
Since $W$ is centered Gaussian, then 
\begin{equation*}
W_{t}-W_{s},
\end{equation*}%
is normally distributed, and its variance is 
\begin{equation*}
\mathbb{E}[(W_{t}-W_{s})^{2}]=\mathbb{E}[W_{t}^{2}]+\mathbb{E}[W_{s}^{2}]-2%
\mathbb{E}[W_{s}W_{t}]=t+s-2s=t-s.
\end{equation*}%
Therefore, 
\begin{equation*}
\ W_{t}-W_{s}\sim N(0,t-s).
\end{equation*}%
Now consider two nonoverlapping increments, where 
\begin{equation*}
0\leq r\leq s\leq t\leq u\leq 1.
\end{equation*}%
Their covariance is 
\begin{align*}
& \func{Cov}(W_{s}-W_{r},W_{u}-W_{t})=\mathbb{E}[(W_{s}-W_{r})(W_{u}-W_{t})]
\\
& =\min (s,u)-\min (s,t)-\min (r,u)+\min (r,t)=s-s-r+r=0.
\end{align*}%
Since $W$ is Gaussian, zero covariance implies independence, and therefore, $%
W$ has independent increments. Combining this with the continuity of its
sample paths, $W_{0}=0$, and 
\begin{equation*}
W_{t}-W_{s}\sim N(0,t-s),
\end{equation*}%
we conclude that $W$ is standard Brownian motion. This completes the proof
of the theorem.
\end{proof}

We have constructed standard Brownian motion using finite stochastic
multivariate Taylor measure functions. The resulting representation is 
\begin{equation*}
W_{t}=Z_{\ast }t+\sum_{j=0}^{\infty
}\sum_{k=0}^{2^{j}-1}Z_{j,k}2^{-(j+2)/2}h_{j,k}(t),
\end{equation*}%
where the finite partial sums are stochastic MTMFs, and their convergence to
Brownian motion is almost surely uniform. Thus, Brownian motion belongs to
the uniform limit completion of the finite stochastic MTMF class.

As defined above, $W_{t}$ is also a martingale.

\begin{corollary}
Let 
\begin{equation*}
\mathbb{F}_{t}^{W}=\sigma (W_{s}:0\leq s\leq t),
\end{equation*}%
be the natural filtration of the constructed process, with its usual
augmentation.

Then 
\begin{equation*}
(W_{t})_{0\leq t\leq 1},
\end{equation*}%
is a martingale and 
\begin{equation*}
(W_{t}^{2}-t)_{0\leq t\leq 1},
\end{equation*}%
is also a martingale.
\end{corollary}

\begin{proof}
For $s\leq t$, independent increments give 
\begin{equation*}
W_{t}=W_{s}+(W_{t}-W_{s}),
\end{equation*}%
where $W_{t}-W_{s}$ is independent of $\mathbb{F}_{s}^{W}$ and has mean
zero. Hence, we have 
\begin{equation*}
\mathbb{E}[W_{t}\mid \mathbb{F}_{s}^{W}]=W_{s}+\mathbb{E}[W_{t}-W_{s}\mid 
\mathbb{F}_{s}^{W}]=W_{s},
\end{equation*}%
and thus, $W$ is a martingale. Furthermore, we can write%
\begin{equation*}
W_{t}^{2}=W_{s}^{2}+2W_{s}(W_{t}-W_{s})+(W_{t}-W_{s})^{2}.
\end{equation*}%
Taking conditional expectations leads to 
\begin{equation*}
\mathbb{E}[W_{t}^{2}\mid \mathbb{F}_{s}^{W}]=W_{s}^{2}+(t-s).
\end{equation*}%
and therefore 
\begin{equation*}
\mathbb{E}[W_{t}^{2}-t\mid \mathbb{F}_{s}^{W}]=W_{s}^{2}-s.
\end{equation*}%
This establishes the second martingale property.
\end{proof}

There are three particularly useful results here.

1. Brownian motion as an MTMF limit: we obtained a concrete example of a
stochastic process generated by finite Taylor measure functions.

2. Uniform convergence: The convergence is stronger than pointwise
convergence. It guarantees continuous sample paths of the limiting process.

3. Gaussian covariance structure: The MTMF construction recovers the
defining Brownian covariance function: 
\begin{equation*}
K(s,t)=\min (s,t).
\end{equation*}%
This construction provides a concrete example of how the stochastic MTMF
framework can be used to represent and approximate fundamental stochastic
processes. A natural question arises: can we represent other important
stochastic processes as limits of SMTMFs, as in the case of Brownian motion?
The answer is yes, and we establish an important characterization of the
SMTMF in the next section.

\section{Completion of the SMTMF Space}

We introduce a uniform limit completion of the stochastic multivariate
Taylor measure function space. The purpose is to establish a function space
that contains the original finite stochastic multivariate Taylor measure
functions and is closed under almost sure uniform convergence.

To that end, we define an appropriate norm, identify the corresponding
equivalence classes of stochastic functions, and establish completeness of
the resulting space. We then prove that fundamental stochastic processes,
such as Brownian motion, Brownian bridge and Fractional Brownian motion, are
constructed through finite stochastic multivariate Taylor measure functions
and belong to this completed space.

This provides a rigorous framework for treating Brownian motion and other
continuous stochastic processes as elements of a generalized stochastic
multivariate Taylor measure function space.

Let $(\Omega ,\mathcal{A},\mathbb{P})$ be a complete probability space, and
let $D=[0,1]^{p},$ for some $p\in \mathbb{N}^{+}$. Let $\mathcal{B}(D)$
denote the Borel $\sigma $-algebra on $D$.

Suppose that $\mathscr S$ denotes the collection of stochastic multivariate
Taylor measure functions of the form 
\begin{equation*}
F(\omega,x) = \sum_{n\in B} a_n(\omega,x)\frac{g(\omega,x)^n}{n!},
\end{equation*}
where $B\in\mathcal{B}(\mathbb{N}_0)$ is finite and the coefficient
functions and Taylor parameter satisfy the appropriate measurability and
finiteness assumptions.

We assume that every element of $\mathscr S$ has a jointly measurable
representative 
\begin{equation*}
F:\Omega \times D\longrightarrow \mathbb{R},
\end{equation*}%
with continuous sample paths in the spatial variable $x$.

The last assumption is imposed because we wish to use the uniform topology
on $C(D)$.

\begin{remark}
If the original definition of $\mathscr S$ does not require continuous
sample paths, one can instead work with bounded measurable functions and an
essential-supremum norm. The present construction is intended for the
continuous-path subclass containing the Brownian approximations.
\end{remark}

\subsection{Construction of the Completed Space}

First, we will work with the essential uniform norm.

\begin{definition}[Essential Uniform Norm]
For a jointly measurable stochastic function 
\begin{equation*}
F:\Omega \times D\longrightarrow \mathbb{R},
\end{equation*}%
with almost surely bounded sample paths, define 
\begin{equation*}
{\Vert F\Vert _{\infty ,\mathrm{ess}}=\limfunc{ess}_{\omega \in \Omega
}\sup_{x\in D}|F(\omega ,x)|.}
\end{equation*}%
Here, the essential supremum is taken with respect to $\mathbb{P}$.
\end{definition}

The quantity above may be infinite for some stochastic functions. We
therefore introduce the bounded subclass 
\begin{equation*}
\mathscr S_{b}=\left\{ F\in \mathscr S:\Vert F\Vert _{\infty ,\mathrm{ess}%
}<\infty \right\} .
\end{equation*}%
We identify two elements $F,G\in \mathscr S_{b}$ whenever 
\begin{equation*}
\Vert F-G\Vert _{\infty ,\mathrm{ess}}=0.
\end{equation*}%
This identification is important because stochastic functions are naturally
defined up to almost sure equality.

\begin{definition}[Uniform-Limit Completion]
Let $\overline{\mathscr S}$ denote the completion of $\mathscr S_{b}$ with
respect to the metric 
\begin{equation*}
d(F,G)=\Vert F-G\Vert _{\infty ,\mathrm{ess}}.
\end{equation*}%
Equivalently, $\overline{\mathscr S}$ consists of equivalence classes of
Cauchy sequences 
\begin{equation*}
(F_{n})_{n\geq 1}\subset \mathscr S_{b},
\end{equation*}%
under the equivalence relation 
\begin{equation*}
(F_{n})\sim (G_{n}),
\end{equation*}%
if 
\begin{equation*}
\lim_{n\rightarrow \infty }\Vert F_{n}-G_{n}\Vert _{\infty ,\mathrm{ess}}=0.
\end{equation*}
\end{definition}

We now establish the completeness result.

\begin{theorem}[Completeness of the Uniform Limit SMTMF Space]
Suppose that every element of $\mathscr S_{b}$ has a jointly measurable
representative with almost surely continuous sample paths on $D=[0,1]^{p}$.
Then the completion $\overline{\mathscr S}$ is a Banach space.

Moreover, every element of $\overline{\mathscr S}$ admits a representative
that is jointly measurable and has almost surely continuous sample paths.
\end{theorem}

\begin{proof}
We divide the proof into several steps.

\noindent \textbf{Step 1. The metric is well-defined}

By definition, 
\begin{equation*}
d(F,G)=\limfunc{ess}_{\omega }\sup_{x\in D}|F(\omega ,x)-G(\omega ,x)|.
\end{equation*}%
Clearly, 
\begin{equation*}
d(F,G)\geq 0,
\end{equation*}%
and 
\begin{equation*}
d(F,G)=d(G,F).
\end{equation*}%
For any $F,G,H\in \mathscr S_{b}$, the triangle inequality gives 
\begin{equation*}
\sup_{x\in D}\left\vert F(\omega ,x)-H(\omega ,x)\right\vert \leq \sup_{x\in
D}\left\vert F(\omega ,x)-G(\omega ,x)\right\vert +\sup_{x\in D}\left\vert
G(\omega ,x)-H(\omega ,x)\right\vert .
\end{equation*}%
Taking essential suprema, we obtain 
\begin{equation*}
d(F,H)\leq d(F,G)+d(G,H).
\end{equation*}%
Finally, after identifying functions that are equal almost surely uniformly
on $D$, 
\begin{equation*}
d(F,G)=0,
\end{equation*}%
implies that $F$ and $G$ represent the same equivalence class. Thus, $d$ is
a metric on the quotient space.

\noindent \textbf{Step 2. Consider a Cauchy sequence}

Let 
\begin{equation*}
(F_{n})_{n\geq 1},
\end{equation*}%
be a Cauchy sequence in $\mathscr S_{b}$. By the definition of a Cauchy
sequence, for every $\varepsilon >0$ there exists $N_{\varepsilon }$ such
that 
\begin{equation*}
d(F_{n},F_{m})<\varepsilon ,
\end{equation*}%
whenever $n,m\geq N_{\varepsilon }$. Choose a subsequence 
\begin{equation*}
(F_{n_{k}})_{k\geq 1},
\end{equation*}%
such that 
\begin{equation*}
d(F_{n_{k+1}},F_{n_{k}})<2^{-k}.
\end{equation*}%
For every $k$, by the definition of the essential supremum, there exists a
null set $A_{k}$ such that 
\begin{equation*}
\sup_{x\in D}\left\vert F_{n_{k+1}}(\omega ,x)-F_{n_{k}}(\omega
,x)\right\vert \leq 2^{-k},
\end{equation*}%
for all $\omega \notin A_{k}$. Let 
\begin{equation*}
A=\bigcup_{k=1}^{\infty }A_{k}.
\end{equation*}%
Since this is a countable union of null sets, we have 
\begin{equation*}
\mathbb{P}(A)=0.
\end{equation*}%
For every $\omega \notin A$, 
\begin{equation*}
\sup_{x\in D}\left\vert F_{n_{k+1}}(\omega ,x)-F_{n_{k}}(\omega
,x)\right\vert \leq 2^{-k}.
\end{equation*}%
Consequently, we have 
\begin{equation*}
\sup_{x\in D}\left\vert F_{n_{\ell }}(\omega ,x)-F_{n_{k}}(\omega
,x)\right\vert \leq \sum_{j=k}^{\ell -1}2^{-j}.
\end{equation*}%
The right-hand side tends to zero as $k,\ell \rightarrow \infty $, and
therefore, for every $\omega \notin A$, the sequence 
\begin{equation*}
(F_{n_{k}}(\omega ,\cdot ))_{k\geq 1},
\end{equation*}%
is uniformly Cauchy on $D$.

\noindent \textbf{Step 3. Define the limiting function}

Since $C(D)$ is complete under the uniform norm, there exists 
\begin{equation*}
F(\omega ,\cdot )\in C(D),
\end{equation*}%
such that 
\begin{equation*}
\lim_{k\rightarrow \infty }\sup_{x\in D}|F_{n_{k}}(\omega ,x)-F(\omega ,x)=0,
\end{equation*}%
for every $\omega \notin A$. Define $F=0$ on $A$. Thus, $F$ has continuous
sample paths almost surely.

\noindent \textbf{Step 4. Establish joint measurability}

The sequence 
\begin{equation*}
F_{n_{k}}:\Omega \times D\rightarrow \mathbb{R},
\end{equation*}%
consists of jointly measurable functions.

Define 
\begin{equation*}
F(\omega ,x)=\limsup_{k\rightarrow \infty }F_{n_{k}}(\omega ,x),
\end{equation*}%
outside $A$, with $F=0$ on $A$. Since the $\lim \sup $ of measurable
functions is measurable, $F$ is jointly measurable. Moreover, the
construction in Step 3 shows that this pointwise limit agrees with the
uniform limit almost surely.

\noindent \textbf{Step 5. Show convergence in the metric}

For almost every $\omega $, we can write%
\begin{equation*}
\sup_{x\in D}\left\vert F_{n_{k}}(\omega ,x)-F(\omega ,x)\right\vert \leq
\sum_{j=k}^{\infty }2^{-j}.
\end{equation*}%
Taking essential suprema gives 
\begin{equation*}
d(F_{n_{k}},F)\leq \sum_{j=k}^{\infty }2^{-j},
\end{equation*}%
and therefore 
\begin{equation*}
d(F_{n_{k}},F)\longrightarrow 0.
\end{equation*}%
Since the original sequence $(F_{n})$ is Cauchy and $(F_{n_{k}})$ converges
to $F$, the entire sequence converges to $F$. Thus, we arrive at%
\begin{equation*}
d(F_{n},F)\longrightarrow 0.
\end{equation*}%
We have shown that every Cauchy sequence has a limit in the completed space,
and as a result, $\overline{\mathscr S}$ is complete.
\end{proof}

\subsection{Application to Brownian Motion}

We now apply the completion theorem to the Brownian motion construction. Let 
$Z_{\ast },Z_{j,k}$ be mutually independent standard normal random variables
and $h_{j,k}$ denote the Schauder functions on $[0,1]$. Define 
\begin{equation*}
F_{N}(t)=Z_{\ast
}t+\sum_{j=0}^{N}\sum_{k=0}^{2^{j}-1}Z_{j,k}2^{-(j+2)/2}h_{j,k}(t),
\end{equation*}%
such that each $F_{N}$ is a finite stochastic MTMF.

We first establish the stronger integrability condition needed for
membership in the essential-uniform completion.

\begin{lemma}[Finite Expected Uniform Norm]
For the Brownian MTMF approximations, 
\begin{equation*}
\mathbb{E}[\Vert F_{N}\Vert _{\infty }]<\infty ,
\end{equation*}%
for every $N\in \mathbb{N}_{0}$. Furthermore, 
\begin{equation*}
\mathbb{E}\left[ \sum_{j=0}^{\infty }2^{-(j+2)/2}\max_{0\leq
k<2^{j}}|Z_{j,k}|\right] <\infty .
\end{equation*}
\end{lemma}

\begin{proof}
For 
\begin{equation*}
M_{j}=\max_{0\leq k<2^{j}}|Z_{j,k}|,
\end{equation*}%
the Gaussian tail estimate gives 
\begin{equation*}
\mathbb{P}(M_{j}>u)\leq \min \{1,2^{j+1}e^{-u^{2}/2}\}.
\end{equation*}%
Using the identity 
\begin{equation*}
\mathbb{E}[M_{j}]=\int_{0}^{\infty }\mathbb{P}(M_{j}>u)\,du,
\end{equation*}%
we obtain 
\begin{equation*}
\mathbb{E}[M_{j}]\leq C_{0}\sqrt{j+1},
\end{equation*}%
for some deterministic constant $C_{0}$. Hence, we can write%
\begin{equation*}
\mathbb{E}\left[ \sum_{j=0}^{\infty }c_{j}M_{j}\right] =\sum_{j=0}^{\infty
}c_{j}\mathbb{E}[M_{j}]\leq C_{0}\sum_{j=0}^{\infty }2^{-(j+2)/2}\sqrt{j+1}%
<\infty .
\end{equation*}%
The finite approximation satisfies 
\begin{equation*}
\Vert F_{N}\Vert _{\infty }\leq |Z_{\ast }|+\sum_{j=0}^{N}c_{j}M_{j}.
\end{equation*}%
Since a standard normal random variable has finite first absolute moment,
the first assertion follows.
\end{proof}

The inclusion theorem for Brownian motion follows.

\begin{theorem}[Brownian Motion Belongs to the Completed Stochastic MTMF
Space]
Let $\mathscr S_{b}$ be the space of bounded finite stochastic MTMFs on $%
[0,1]$, equipped with the essential-uniform metric 
\begin{equation*}
d(F,G)=\limfunc{ess}_{\omega }\sup_{t\in \lbrack 0,1]}|F(\omega ,t)-G(\omega
,t)|.
\end{equation*}%
Let $\overline{\mathscr S}$ be its uniform-limit completion. Suppose that
the finite Brownian MTMF approximations are 
\begin{equation*}
F_{N}(t)=Z_{\ast
}t+\sum_{j=0}^{N}\sum_{k=0}^{2^{j}-1}Z_{j,k}2^{-(j+2)/2}h_{j,k}(t).
\end{equation*}%
Then the sequence $(F_{N})$ converges in the metric $d$ to a continuous
stochastic process $W$. Moreover, $W$ belongs to $\overline{\mathscr S}$ and
is standard Brownian motion.
\end{theorem}

\begin{proof}
For $N>M$, 
\begin{equation*}
\Vert F_{N}-F_{M}\Vert _{\infty }\leq \sum_{j=M+1}^{N}c_{j}M_{j}.
\end{equation*}%
Taking expectations gives 
\begin{equation*}
\mathbb{E}[\Vert F_{N}-F_{M}\Vert _{\infty }]\leq \sum_{j=M+1}^{N}c_{j}%
\mathbb{E}[M_{j}].
\end{equation*}%
By the preceding lemma, 
\begin{equation*}
\sum_{j=0}^{\infty }c_{j}\mathbb{E}[M_{j}]<\infty .
\end{equation*}%
Consequently, 
\begin{equation*}
\lim_{M,N\rightarrow \infty }\mathbb{E}[\Vert F_{N}-F_{M}\Vert _{\infty }]=0.
\end{equation*}%
Thus, $(F_{N})$ is Cauchy in the essential-uniform metric.

Since $\overline{\mathscr S}$ is complete, there exists $W\in \overline{%
\mathscr S}$ such that 
\begin{equation*}
d(F_{N},W)\longrightarrow 0.
\end{equation*}%
The almost sure uniform convergence established in the Brownian construction
identifies this limit with 
\begin{equation*}
W_{t}=Z_{\ast }t+\sum_{j=0}^{\infty
}\sum_{k=0}^{2^{j}-1}Z_{j,k}2^{-(j+2)/2}h_{j,k}(t).
\end{equation*}%
The Brownian construction theorem establishes that this process has
continuous sample paths, centered Gaussian finite-dimensional distributions,
and covariance 
\begin{equation*}
\mathbb{E}[W_{s}W_{t}]=\min (s,t).
\end{equation*}%
It also has independent increments and starts at zero. Therefore, $W$ is
standard Brownian motion.

Since $W$ is the limit in $\overline{\mathscr S}$ of a sequence of elements
of $\mathscr S_b$, we conclude that 
\begin{equation*}
{W\in\overline{\mathscr S}.}
\end{equation*}
\end{proof}

The completion also preserves important properties of the original
stochastic MTMF functions under suitable convergence conditions.

\begin{proposition}[Preservation of Adaptedness]
Let $(\mathcal{A}_{t})_{t\geq 0}$ be a filtration. Suppose that every $F_{N}$
is adapted, meaning that for each $t$, 
\begin{equation*}
F_{N}(t)\in \mathcal{A}_{t}.
\end{equation*}%
If 
\begin{equation*}
F_{N}(t)\longrightarrow F(t)\quad \text{almost surely,}
\end{equation*}%
for every fixed $t$, then $F(t)$ is $\mathcal{A}_{t}$-measurable whenever $%
\mathcal{A}_{t}$ is complete. Consequently, the limiting process is adapted.
\end{proposition}

\begin{proof}
Fix $t\geq 0$. Each random variable $F_{N}(t)$ is $\mathcal{A}_{t}$%
-measurable. The almost sure limit 
\begin{equation*}
F(t)=\lim_{N\rightarrow \infty }F_{N}(t),
\end{equation*}%
is therefore measurable with respect to the completion of $\mathcal{A}_{t}$.
Since $\mathcal{A}_{t}$ is complete, the limit is $\mathcal{A}_{t}$%
-measurable. This holds for every fixed $t$, so $F$ is adapted.
\end{proof}

\begin{remark}[A Filtration Issue for Brownian Motion]
In the Brownian construction, the Gaussian coefficients $Z_{j,k}$ need not
individually be measurable with respect to the Brownian natural filtration
at time $t$.

Nevertheless, the process $W_t$ is adapted to its own natural filtration by
definition.

Thus, adaptedness of the limit should be verified using the limiting process
itself when the approximating coefficients are defined on a larger
probability space.
\end{remark}

\subsection{A Stronger Integrability Based Alternative}

The essential uniform norm is particularly strong. It requires an almost
surely deterministic uniform bound and is not appropriate for every
stochastic process. An alternative is to use the norm 
\begin{equation*}
\Vert F\Vert _{1,\infty }=E\left[ \sup_{x\in D}|F(\omega ,x)|\right] .
\end{equation*}%
Define 
\begin{equation*}
\mathscr S_{1}=\left\{ F\in \mathscr S:\mathbb{E}[\Vert F\Vert _{\infty
}]<\infty \right\} .
\end{equation*}%
The corresponding completion is denoted by $\overline{\mathscr S}_{1}$.

\begin{theorem}[Brownian Membership Under the Expected Uniform Norm]
Let $\overline{\mathscr S}_{1}$ denote the completion of $\mathscr S_{1}$
under the metric 
\begin{equation*}
d_{1}(F,G)=\mathbb{E}\left[ \sup_{t\in \lbrack 0,1]}|F(t)-G(t)|\right] .
\end{equation*}%
Then the Brownian MTMF approximations 
\begin{equation*}
F_{N}(t)=Z_{\ast
}t+\sum_{j=0}^{N}\sum_{k=0}^{2^{j}-1}Z_{j,k}2^{-(j+2)/2}h_{j,k}(t),
\end{equation*}%
converge in $\overline{\mathscr S}_{1}$ to standard Brownian motion.
\end{theorem}

\begin{proof}
For $N>M$, we have%
\begin{equation*}
d_{1}(F_{N},F_{M})=\mathbb{E}\left[ \sup_{t\in \lbrack
0,1]}|F_{N}(t)-F_{M}(t)|\right] \leq \sum_{j=M+1}^{N}2^{-(j+2)/2}\mathbb{E}%
[M_{j}].
\end{equation*}%
As established previously, 
\begin{equation*}
\mathbb{E}[M_{j}]\leq C_{0}\sqrt{j+1}.
\end{equation*}%
Therefore, we can write%
\begin{equation*}
d_{1}(F_{N},F_{M})\leq C_{0}\sum_{j=M+1}^{N}2^{-(j+2)/2}\sqrt{j+1}.
\end{equation*}%
Since 
\begin{equation*}
\sum_{j=0}^{\infty }2^{-(j+2)/2}\sqrt{j+1}<\infty ,
\end{equation*}%
the right-hand side tends to zero as $M,N\rightarrow \infty $. Thus, $%
(F_{N}) $ is Cauchy in the $d_{1}$ metric. By completeness of the
corresponding function space, there exists a limit $W\in \overline{\mathscr S%
}_{1}$. The almost sure uniform convergence identifies this limit with the
Brownian motion constructed abov, and hence 
\begin{equation*}
{W\in \overline{\mathscr S}_{1}.}
\end{equation*}
\end{proof}

We have introduced a completion of the stochastic multivariate Taylor
measure function space that allows limits of finite stochastic MTMFs to
remain within a well defined function space. In particular, we have shown
that standard Brownian motion belongs to the completed space: 
\begin{equation*}
{W\in \overline{\mathscr S}.}
\end{equation*}%
We have also established membership under the weaker expected-uniform norm: 
\begin{equation*}
{W\in \overline{\mathscr S}_{1}.}
\end{equation*}%
The completion preserves continuity and joint measurability under the stated
assumptions, and adaptedness under the additional filtration assumptions.
This provides a rigorous mathematical framework for extending the stochastic
MTMF theory beyond finite Taylor measures to include continuous time
stochastic processes such as Brownian motion.

One further point for mathematical rigor; the completion should formally be
taken after identifying processes equal almost surely as elements of $%
C([0,1])$. This makes the resulting space a genuine Banach space and avoids
ambiguity about different versions of the same stochastic process.

\subsection{Application to Brownian Bridge}

We establish that a standard Brownian bridge can be represented as a limit
of finite stochastic multivariate Taylor measure functions (MTMFs). In
particular, we construct a sequence of finite stochastic MTMFs using
piecewise linear interpolation and establish almost sure uniform convergence
and convergence in the expected uniform norm. Consequently, the Brownian
bridge belongs to the completion of the finite stochastic MTMF space.

As before, we take $(\Omega ,\mathcal{A},\mathbb{P})$ to be a complete
probability space supporting a standard Brownian motion 
\begin{equation*}
W=\{W_{t}:t\in \lbrack 0,T]\},
\end{equation*}%
where $T>0$. Recall that a finite stochastic MTMF has the representation 
\begin{equation*}
F(\omega ,t)=\sum_{k\in B}a_{k}(\omega ,t)\frac{g(\omega ,t)^{k}}{k!},
\end{equation*}%
where $B\subseteq \mathbb{N}_{0}$ is finite, and the coefficient functions
and Taylor parameter satisfy the required measurability and finiteness
conditions.

We assume that the MTMF space permits jointly measurable random coefficient
functions of $(\omega ,t)$. Let $\mathscr S$ denote the space of continuous
path finite stochastic MTMFs.

Define the expected-uniform norm by 
\begin{equation*}
\|F\|_{1,\infty} = \mathbb{E}\left[ \sup_{t\in[0,T]}|F(\omega,t)| \right].
\end{equation*}

We denote by 
\begin{equation*}
\overline{\mathscr S}_{1},
\end{equation*}%
the completion of $\mathscr S$ with respect to this norm.

For each $n\geq1$, consider the dyadic partition 
\begin{equation*}
t_{n,k}=\frac{kT}{2^n}, \qquad k=0,\ldots,2^n.
\end{equation*}

Let $\varphi _{n,k}$ denote the associated piecewise-linear hat functions.
These functions satisfy 
\begin{equation*}
\varphi _{n,k}(t)\geq 0,
\end{equation*}%
\begin{equation*}
\sum_{k=0}^{2^{n}}\varphi _{n,k}(t)=1,
\end{equation*}%
and 
\begin{equation*}
\varphi _{n,k}(t)\neq 0\implies \left\vert t-t_{n,k}\right\vert \leq \frac{T%
}{2^{n}}.
\end{equation*}%
Now recall the definition of the Brownian bridge.

\begin{definition}
The standard Brownian bridge on $[0,T]$ is the stochastic process 
\begin{equation*}
B_{t}^{\circ }=W_{t}-\frac{t}{T}W_{T},\qquad 0\leq t\leq T.
\end{equation*}%
It is a centered Gaussian process with covariance function 
\begin{equation*}
R^{\circ }(s,t)=\min (s,t)-\frac{st}{T}.
\end{equation*}%
Moreover, 
\begin{equation*}
B_{0}^{\circ }=B_{T}^{\circ }=0\quad \text{almost surely}.
\end{equation*}
\end{definition}

The inclusion theorem follows.

\begin{theorem}[Brownian Bridge Representation in the Completed MTMF Space]
Let $B^{\circ }$ be a standard Brownian bridge on $[0,T]$. For every $n\geq
1 $, define 
\begin{equation*}
F_{n}^{\circ }(\omega ,t)=\sum_{k=0}^{2^{n}}B_{t_{n,k}}^{\circ }(\omega
)\varphi _{n,k}(t).
\end{equation*}%
Then the following statements hold.

\begin{enumerate}
\item Every $F_{n}^{\circ }$ is a finite stochastic MTMF.

\item Every $F_{n}^{\circ }$ is jointly measurable on $\Omega \times \lbrack
0,T]$.

\item The approximating sequence converges almost surely uniformly: 
\begin{equation*}
\sup_{t\in \lbrack 0,T]}|F_{n}^{\circ }(t)-B_{t}^{\circ }|\longrightarrow
0\quad \text{almost surely}.
\end{equation*}

\item The convergence also holds in the expected-uniform norm: 
\begin{equation*}
\lim_{n\rightarrow \infty }\mathbb{E}\left[ \sup_{t\in \lbrack
0,T]}|F_{n}^{\circ }(t)-B_{t}^{\circ }|\right] =0.
\end{equation*}

Consequently, we have 
\begin{equation*}
{B^{\circ }\in \overline{\mathscr S}_{1}.}
\end{equation*}
\end{enumerate}
\end{theorem}

\begin{proof}
We prove the four assertions separately.

\noindent \textbf{Step 1. Express the approximating process as a finite
stochastic MTMF.}

Recall that 
\begin{equation*}
B_{t_{n,k}}^{\circ }=W_{t_{n,k}}-\frac{t_{n,k}}{T}W_{T}.
\end{equation*}%
Therefore, 
\begin{equation*}
F_{n}^{\circ }(\omega ,t)=\sum_{k=0}^{2^{n}}\left( W_{t_{n,k}}(\omega )-%
\frac{t_{n,k}}{T}W_{T}(\omega )\right) \varphi _{n,k}(t).
\end{equation*}%
Choose the deterministic Taylor parameter 
\begin{equation*}
g(\omega ,t)=1.
\end{equation*}%
For every $k=0,\ldots ,2^{n}$, define 
\begin{equation*}
a_{k}(\omega ,t)=k!\left( W_{t_{n,k}}(\omega )-\frac{t_{n,k}}{T}W_{T}(\omega
)\right) \varphi _{n,k}(t).
\end{equation*}%
Then 
\begin{align*}
\sum_{k=0}^{2^{n}}a_{k}(\omega ,t)\frac{g(\omega ,t)^{k}}{k!}&
=\sum_{k=0}^{2^{n}}k!B_{t_{n,k}}^{\circ }(\omega )\varphi _{n,k}(t)\frac{1}{%
k!} \\
& =\sum_{k=0}^{2^{n}}B_{t_{n,k}}^{\circ }(\omega )\varphi
_{n,k}(t)=F_{n}^{\circ }(\omega ,t).
\end{align*}%
Since the summation contains only finitely many terms, $F_{n}^{\circ }$ is a
finite stochastic MTMF, and thus 
\begin{equation*}
F_{n}^{\circ }\in \mathscr S.
\end{equation*}

\noindent \textbf{Step 2. Establish joint measurability.}

For each fixed $k$, the mapping 
\begin{equation*}
\omega \longmapsto B_{t_{n,k}}^{\circ }(\omega ),
\end{equation*}%
is $\mathcal{A}$-measurable because $B_{t_{n,k}}^{\circ }$ is a random
variable. Furthermore, 
\begin{equation*}
t\longmapsto \varphi _{n,k}(t),
\end{equation*}%
is deterministic and continuous. Consequently, the mapping 
\begin{equation*}
(\omega ,t)\longmapsto B_{t_{n,k}}^{\circ }(\omega )\varphi _{n,k}(t),
\end{equation*}%
is measurable with respect to 
\begin{equation*}
\mathcal{A}\otimes \mathcal{B}([0,T]).
\end{equation*}%
Since $F_{n}^{\circ }$ is a finite sum of jointly measurable functions, we
conclude that 
\begin{equation*}
F_{n}^{\circ }:\Omega \times \lbrack 0,T]\longrightarrow \mathbb{R}
\end{equation*}%
is jointly measurable.

\noindent \textbf{Step 3. Prove almost sure uniform convergence.}

For a continuous function $f:[0,T]\rightarrow \mathbb{R}$, define its
modulus of continuity by 
\begin{equation*}
\omega _{f}(\delta )=\sup_{\substack{ s,t\in \lbrack 0,T]  \\ |s-t|\leq
\delta }}|f(s)-f(t)|.
\end{equation*}%
Since Brownian motion has almost surely continuous sample paths and 
\begin{equation*}
B_{t}^{\circ }=W_{t}-\frac{t}{T}W_{T},
\end{equation*}%
the Brownian bridge also has almost surely continuous sample paths. Fix $%
\omega $ outside the null set on which continuity fails. We have 
\begin{equation*}
F_{n}^{\circ }(\omega ,t)-B_{t}^{\circ }(\omega
)=\sum_{k=0}^{2^{n}}B_{t_{n,k}}^{\circ }(\omega )\varphi
_{n,k}(t)-B_{t}^{\circ }(\omega )=\sum_{k=0}^{2^{n}}\varphi _{n,k}(t)\left(
B_{t_{n,k}}^{\circ }(\omega )-B_{t}^{\circ }(\omega )\right) ,
\end{equation*}%
where we used 
\begin{equation*}
\sum_{k=0}^{2^{n}}\varphi _{n,k}(t)=1.
\end{equation*}%
Taking absolute values and applying the triangle inequality gives 
\begin{equation*}
\left\vert F_{n}^{\circ }(\omega ,t)-B_{t}^{\circ }(\omega )\right\vert \leq
\sum_{k=0}^{2^{n}}\varphi _{n,k}(t)|B_{t_{n,k}}^{\circ }(\omega
)-B_{t}^{\circ }(\omega )|.
\end{equation*}%
Whenever 
\begin{equation*}
\varphi _{n,k}(t)\neq 0,
\end{equation*}%
we have 
\begin{equation*}
\left\vert t-t_{n,k}\right\vert \leq \frac{T}{2^{n}}.
\end{equation*}%
Therefore, 
\begin{equation*}
\left\vert F_{n}^{\circ }(\omega ,t)-B_{t}^{\circ }(\omega )\right\vert \leq
\omega _{B^{\circ }(\omega )}\left( \frac{T}{2^{n}}\right)
\sum_{k=0}^{2^{n}}\varphi _{n,k}(t)=\omega _{B^{\circ }(\omega )}\left( 
\frac{T}{2^{n}}\right) .
\end{equation*}%
Taking the supremum over $t\in \lbrack 0,T]$, we obtain 
\begin{equation*}
\sup_{t\in \lbrack 0,T]}\left\vert F_{n}^{\circ }(\omega ,t)-B_{t}^{\circ
}(\omega )\right\vert \leq \omega _{B^{\circ }(\omega )}\left( \frac{T}{2^{n}%
}\right) .
\end{equation*}%
Because every continuous function on a compact interval is uniformly
continuous, we have%
\begin{equation*}
\lim_{n\rightarrow \infty }\omega _{B^{\circ }(\omega )}\left( \frac{T}{2^{n}%
}\right) =0,
\end{equation*}%
and hence 
\begin{equation*}
\sup_{t\in \lbrack 0,T]}|F_{n}^{\circ }(t)-B_{t}^{\circ }|\longrightarrow
0\quad \text{almost surely}.
\end{equation*}

\noindent \textbf{Step 4. Prove convergence in the expected-uniform norm.}

We first establish an integrable dominating random variable. By definition, 
\begin{equation*}
B_{t}^{\circ }=W_{t}-\frac{t}{T}W_{T}.
\end{equation*}%
Consequently, 
\begin{equation*}
\sup_{t\in \lbrack 0,T]}\left\vert B_{t}^{\circ }\right\vert \leq \sup_{t\in
\lbrack 0,T]}|W_{t}|+\sup_{t\in \lbrack 0,T]}\frac{t}{T}\left\vert
W_{T}\right\vert \leq \sup_{t\in \lbrack 0,T]}|W_{t}|+\left\vert
W_{T}\right\vert \leq 2\sup_{t\in \lbrack 0,T]}|W_{t}|.
\end{equation*}%
The Brownian supremum has finite expectation. Indeed, the reflection
principle gives 
\begin{equation*}
\mathbb{E}\left[ \sup_{t\in \lbrack 0,T]}|W_{t}|\right] <\infty .
\end{equation*}%
Therefore, 
\begin{equation*}
\mathbb{E}\left[ \sup_{t\in \lbrack 0,T]}|B_{t}^{\circ }|\right] <\infty .
\end{equation*}%
Next, because the interpolation weights are nonnegative and sum to one, 
\begin{equation*}
\left\vert F_{n}^{\circ }(t)\right\vert \leq \sum_{k=0}^{2^{n}}\varphi
_{n,k}(t)\left\vert B_{t_{n,k}}^{\circ }\right\vert \leq \sup_{s\in \lbrack
0,T]}|B_{s}^{\circ }|.
\end{equation*}%
It follows that 
\begin{equation*}
\sup_{t\in \lbrack 0,T]}\left\vert F_{n}^{\circ }(t)-B_{t}^{\circ
}\right\vert \leq \sup_{t\in \lbrack 0,T]}\left\vert F_{n}^{\circ
}(t)\right\vert +\sup_{t\in \lbrack 0,T]}\left\vert B_{t}^{\circ
}\right\vert \leq 2\sup_{t\in \lbrack 0,T]}|B_{t}^{\circ }|.
\end{equation*}%
The right-hand side is integrable and does not depend on $n$. By Step 3, the
random variables 
\begin{equation*}
\sup_{t\in \lbrack 0,T]}|F_{n}^{\circ }(t)-B_{t}^{\circ }|,
\end{equation*}%
converge almost surely to zero. The dominated convergence theorem therefore
implies 
\begin{equation*}
\lim_{n\rightarrow \infty }E\left[ \sup_{t\in \lbrack 0,T]}|F_{n}^{\circ
}(t)-B_{t}^{\circ }|\right] =0.
\end{equation*}

\noindent \textbf{Step 5. Identify the covariance function of the limiting
process.}

Since 
\begin{equation*}
B_{t}^{\circ }=W_{t}-\frac{t}{T}W_{T},
\end{equation*}%
we calculate 
\begin{equation*}
\mathbb{E}[B_{s}^{\circ }B_{t}^{\circ }]=\mathbb{E}[W_{s}W_{t}]-\frac{t}{T}%
\mathbb{E}[W_{s}W_{T}]-\frac{s}{T}\mathbb{E}[W_{t}W_{T}]+\frac{st}{T^{2}}%
\mathbb{E}[W_{T}^{2}].
\end{equation*}%
Using the Brownian covariance identities 
\begin{equation*}
\mathbb{E}[W_{s}W_{t}]=\min (s,t),
\end{equation*}%
\begin{equation*}
\mathbb{E}[W_{s}W_{T}]=s,\qquad \mathbb{E}[W_{t}W_{T}]=t,
\end{equation*}%
and 
\begin{equation*}
\mathbb{E}[W_{T}^{2}]=T,
\end{equation*}%
we obtain 
\begin{equation*}
\mathbb{E}[B_{s}^{\circ }B_{t}^{\circ }]=\min (s,t)-\frac{st}{T}-\frac{st}{T}%
+\frac{st}{T}=\min (s,t)-\frac{st}{T}.
\end{equation*}%
Thus, 
\begin{equation*}
\func{Cov}(B_{s}^{\circ },B_{t}^{\circ })=\min (s,t)-\frac{st}{T}.
\end{equation*}

This confirms that the limiting process is the standard Brownian bridge.

\noindent \textbf{Step 6. Conclude membership in the completed stochastic
MTMF space.}

We have established that 
\begin{equation*}
F_{n}^{\circ }\in \mathscr S,
\end{equation*}%
for every $n\geq 1$ and that 
\begin{equation*}
\lim_{n\rightarrow \infty }\Vert F_{n}^{\circ }-B^{\circ }\Vert _{1,\infty
}=0.
\end{equation*}%
Therefore, by the definition of the completion, 
\begin{equation*}
B^{\circ }\in \overline{\mathscr S}_{1}.
\end{equation*}%
This completes the proof.
\end{proof}

The Brownian bridge admits an explicit approximation by finite stochastic
MTMFs: 
\begin{equation*}
F_{n}^{\circ }(t)=\sum_{k=0}^{2^{n}}\left( W_{t_{n,k}}-\frac{t_{n,k}}{T}%
W_{T}\right) \varphi _{n,k}(t).
\end{equation*}%
The approximation converges almost surely uniformly and in the expected
uniform norm. Consequently, the Brownian bridge is an element of the
completed stochastic MTMF space 
\begin{equation*}
B^{\circ }\in \overline{\mathscr S}_{1}.
\end{equation*}

\begin{remark}
The conclusion concerns the completion $\overline{\mathscr S}_1$, not
necessarily the original finite MTMF space $\mathscr S$. This distinction is
important because each approximating function has finitely many Taylor
terms, whereas the Brownian bridge is obtained as a limit of these functions.
\end{remark}

\subsection{Application to Fractional Brownian Motion}

Fractional Brownian motion can be studied through its Hurst parameter and
covariance structure. We establish that fractional Brownian motion with
arbitrary Hurst parameter $H\in (0,1)$ can be represented as a limit of
finite SMTMFs. We construct an explicit sequence of finite stochastic MTMFs
using piecewise-linear interpolation over dyadic partitions, and prove that
the resulting sequence converges almost surely uniformly and in the expected
uniform norm to fractional Brownian motion. Consequently, fractional
Brownian motion belongs to the completion of the finite stochastic MTMF
space.

Let \ $(\Omega ,\mathcal{A},\mathbb{P})$ be a complete probability space,
and define 
\begin{equation*}
B^{H}=\{B_{t}^{H}:t\in \lbrack 0,T]\},
\end{equation*}%
be a fractional Brownian motion with Hurst parameter 
\begin{equation*}
H\in (0,1).
\end{equation*}%
We assume that the underlying probability space supports a continuous
modification of this process. Recall that a finite stochastic MTMF has the
representation 
\begin{equation*}
F(\omega ,t)=\sum_{k\in A}a_{k}(\omega ,t)\frac{g(\omega ,t)^{k}}{k!},
\end{equation*}%
where $A\subseteq \mathbb{N}_{0}$ is finite, and the coefficient functions
and Taylor parameter satisfy the measurability and finiteness conditions of
the stochastic MTMF space. As before, we assume that the MTMF space permits
jointly measurable random coefficient functions of $(\omega ,t)$.

Let $\mathscr S$ denote the space of continuous-path finite stochastic
MTMFs. Define the expected uniform norm by 
\begin{equation*}
\Vert F\Vert _{1,\infty }=\mathbb{E}\left[ \sup_{t\in \lbrack 0,T]}|F(\omega
,t)|\right] ,
\end{equation*}%
and let 
\begin{equation*}
\overline{\mathscr S}_{1}
\end{equation*}%
denote the completion of $\mathscr S$ under this norm.

For each integer $n\geq 1$, consider the dyadic partition 
\begin{equation*}
t_{n,k}=\frac{kT}{2^{n}},\qquad k=0,\ldots ,2^{n},
\end{equation*}%
and let 
\begin{equation*}
\varphi _{n,k}:[0,T]\longrightarrow \lbrack 0,1],
\end{equation*}%
be the associated piecewise-linear hat functions. These functions satisfy 
\begin{equation*}
\varphi _{n,k}(t)\geq 0,
\end{equation*}%
\begin{equation*}
\sum_{k=0}^{2^{n}}\varphi _{n,k}(t)=1,
\end{equation*}%
and 
\begin{equation*}
\varphi _{n,k}(t)\neq 0\implies \left\vert t-t_{n,k}\right\vert \leq \frac{T%
}{2^{n}}.
\end{equation*}

Recall the definition of fractional Brownian motion.

\begin{definition}
A fractional Brownian motion with Hurst parameter $H\in (0,1)$ is a centered
Gaussian process satisfying 
\begin{equation*}
B_{0}^{H}=0,
\end{equation*}%
and having covariance function 
\begin{equation*}
{\ R_{H}(s,t)=\mathbb{E}[B_{s}^{H}B_{t}^{H}]=\frac{1}{2}\left(
s^{2H}+t^{2H}-|t-s|^{2H}\right) .}
\end{equation*}%
In particular, its increments satisfy 
\begin{equation*}
\mathbb{E}\left[ (B_{t}^{H}-B_{s}^{H})^{2}\right] =|t-s|^{2H}.
\end{equation*}
\end{definition}

The inclusion theorem is given next.

\begin{theorem}[Fractional Brownian Motion inclusion in Completed SMTMF space%
]
Let $B^{H}$ be fractional Brownian motion on $[0,T]$ with Hurst parameter 
\begin{equation*}
H\in (0,1).
\end{equation*}%
For every $n\geq 1$, define 
\begin{equation*}
{\ F_{n}^{H}(\omega ,t)=\sum_{k=0}^{2^{n}}B_{t_{n,k}}^{H}(\omega )\varphi
_{n,k}(t).}
\end{equation*}%
Then the following assertions hold.

\begin{enumerate}
\item Every $F_n^H$ is a finite stochastic MTMF.

\item Every $F_n^H$ is jointly measurable on 
\begin{equation*}
\Omega\times[0,T].
\end{equation*}

\item The approximating sequence converges almost surely uniformly: 
\begin{equation*}
\lim_{n\to\infty} \sup_{t\in[0,T]} |F_n^H(t)-B_t^H|=0 \quad\text{almost
surely}.
\end{equation*}

\item The convergence holds in the expected-uniform norm: 
\begin{equation*}
\lim_{n\to\infty} \mathbb{E}\left[ \sup_{t\in[0,T]} |F_n^H(t)-B_t^H| \right]%
=0.
\end{equation*}

Consequently, 
\begin{equation*}
{B^{H}\in \overline{\mathscr S}_{1}.}
\end{equation*}
\end{enumerate}
\end{theorem}

\begin{proof}
We establish each assertion separately.

\noindent \textbf{Step 1. Construct a finite stochastic MTMF representation}

For each $n\geq 1$, define 
\begin{equation*}
F_{n}^{H}(\omega ,t)=\sum_{k=0}^{2^{n}}B_{t_{n,k}}^{H}(\omega )\varphi
_{n,k}(t).
\end{equation*}%
We choose the deterministic Taylor parameter 
\begin{equation*}
g(\omega ,t)=1.
\end{equation*}%
For every 
\begin{equation*}
k=0,\ldots ,2^{n},
\end{equation*}%
define the random coefficient function 
\begin{equation*}
a_{n,k}(\omega ,t)=k!B_{t_{n,k}}^{H}(\omega )\varphi _{n,k}(t).
\end{equation*}%
Then 
\begin{equation*}
\sum_{k=0}^{2^{n}}a_{n,k}(\omega ,t)\frac{g(\omega ,t)^{k}}{k!}%
=\sum_{k=0}^{2^{n}}k!B_{t_{n,k}}^{H}(\omega )\varphi _{n,k}(t)\frac{1}{k!}%
=\sum_{k=0}^{2^{n}}B_{t_{n,k}}^{H}(\omega )\varphi
_{n,k}(t)=F_{n}^{H}(\omega ,t).
\end{equation*}%
Since the summation contains only finitely many terms, $F_{n}^{H}$ is a
finite stochastic MTMF, and therefore 
\begin{equation*}
F_{n}^{H}\in \mathscr S.
\end{equation*}

\noindent \textbf{Step 2. Establish joint measurability}

For every fixed $k$, the mapping 
\begin{equation*}
\omega \longmapsto B_{t_{n,k}}^{H}(\omega ),
\end{equation*}%
is $\mathcal{A}$-measurable because $B_{t_{n,k}}^{H}$ is a random variable.
Furthermore, the function 
\begin{equation*}
t\longmapsto \varphi _{n,k}(t),
\end{equation*}%
is deterministic and continuous. Consequently, 
\begin{equation*}
(\omega ,t)\longmapsto B_{t_{n,k}}^{H}(\omega )\varphi _{n,k}(t),
\end{equation*}%
is measurable with respect to 
\begin{equation*}
\mathcal{A}\otimes \mathcal{B}([0,T]).
\end{equation*}%
Since $F_{n}^{H}$ is a finite sum of these jointly measurable functions, it
follows that 
\begin{equation*}
F_{n}^{H}:\Omega \times \lbrack 0,T]\longrightarrow \mathbb{R},
\end{equation*}%
is jointly measurable.

\noindent \textbf{Step 3. Establish continuity of fractional Brownian motion}

We recall that fractional Brownian motion has Gaussian increments satisfying 
\begin{equation*}
B_{t}^{H}-B_{s}^{H}\sim N(0,|t-s|^{2H}).
\end{equation*}%
Let $q>0$. Since a centered Gaussian random variable with variance $%
|t-s|^{2H}$ has the same distribution as $|t-s|^{H}Z$, where $Z$ is standard
normal, we obtain 
\begin{equation*}
\mathbb{E}[|B_{t}^{H}-B_{s}^{H}|^{q}]=|t-s|^{qH}\mathbb{E}%
[|Z|^{q}]=C_{q}|t-s|^{qH},
\end{equation*}%
where 
\begin{equation*}
C_{q}=\mathbb{E}[|Z|^{q}]<\infty .
\end{equation*}%
Choose $q$ sufficiently large so that 
\begin{equation*}
qH>1.
\end{equation*}%
For example, any 
\begin{equation*}
q>\frac{1}{H},
\end{equation*}%
is sufficient. Kolmogorov's continuity theorem (see \cite{Kallenberg2002}, 
\cite{AdlerTaylor2007}, \cite{micheas2018theory}) then guarantees the
existence of a continuous modification of $B^{H}$.

Moreover, the modification can be chosen to be locally H\"{o}lder continuous
of every order 
\begin{equation*}
\alpha <H.
\end{equation*}%
Thus, throughout the remainder of the proof, we work with the continuous
modification of fractional Brownian motion.

\noindent \textbf{Step 4. Establish almost sure uniform convergence}

For a continuous function $f:[0,T]\rightarrow \mathbb{R}$, define its
modulus of continuity by 
\begin{equation*}
\omega _{f}(\delta )=\sup_{\substack{ s,t\in \lbrack 0,T]  \\ |s-t|\leq
\delta }}|f(s)-f(t)|.
\end{equation*}%
By construction, 
\begin{equation*}
F_{n}^{H}(t)=\sum_{k=0}^{2^{n}}B_{t_{n,k}}^{H}\varphi _{n,k}(t).
\end{equation*}%
Using the identity 
\begin{equation*}
\sum_{k=0}^{2^{n}}\varphi _{n,k}(t)=1,
\end{equation*}%
we can write 
\begin{equation*}
F_{n}^{H}(t)-B_{t}^{H}=\sum_{k=0}^{2^{n}}\varphi _{n,k}(t)\left(
B_{t_{n,k}}^{H}-B_{t}^{H}\right) .
\end{equation*}%
Taking absolute values gives 
\begin{equation*}
\left\vert F_{n}^{H}(t)-B_{t}^{H}\right\vert \leq \sum_{k=0}^{2^{n}}\varphi
_{n,k}(t)|B_{t_{n,k}}^{H}-B_{t}^{H}|.
\end{equation*}%
Whenever 
\begin{equation*}
\varphi _{n,k}(t)\neq 0,
\end{equation*}%
we have 
\begin{equation*}
\left\vert t-t_{n,k}\right\vert \leq \frac{T}{2^{n}}.
\end{equation*}%
Therefore, 
\begin{equation*}
\left\vert F_{n}^{H}(t)-B_{t}^{H}\right\vert \leq \omega _{B^{H}}\left( 
\frac{T}{2^{n}}\right) \sum_{k=0}^{2^{n}}\varphi _{n,k}(t)=\omega
_{B^{H}}\left( \frac{T}{2^{n}}\right) .
\end{equation*}%
Taking the supremum over $t\in \lbrack 0,T]$ yields 
\begin{equation*}
\Vert F_{n}^{H}-B^{H}\Vert _{\infty }\leq \omega _{B^{H}}\left( \frac{T}{%
2^{n}}\right) .
\end{equation*}%
Since the sample paths of $B^{H}$ are almost surely continuous on the
compact interval $[0,T]$, they are uniformly continuous almost surely.
Consequently, 
\begin{equation*}
\lim_{n\rightarrow \infty }\omega _{B^{H}}\left( \frac{T}{2^{n}}\right)
=0\quad \text{almost surely},
\end{equation*}%
and it follows that 
\begin{equation*}
\lim_{n\rightarrow \infty }\sup_{t\in \lbrack
0,T]}|F_{n}^{H}(t)-B_{t}^{H}|=0\quad \text{almost surely}.
\end{equation*}

\noindent \textbf{Step 5. Establish integrability of the uniform norm}

To obtain convergence in the expected-uniform norm, we first show that 
\begin{equation*}
\mathbb{E}[\Vert B^{H}\Vert _{\infty }]<\infty .
\end{equation*}%
The process $B^{H}$ is a centered Gaussian random element of the separable
Banach space $C([0,T])$, equipped with the uniform norm. Its sample paths
are continuous almost surely. By Fernique's theorem (see \cite{Fernique1970}%
, \cite{Bogachev1998}), there exists a constant $c>0$ such that 
\begin{equation*}
\mathbb{E}\left[ \exp \left( c\Vert B^{H}\Vert _{\infty }^{2}\right) \right]
<\infty .
\end{equation*}%
In particular, 
\begin{equation*}
E\left[ \sup_{t\in \lbrack 0,T]}|B_{t}^{H}|\right] <\infty .
\end{equation*}

\noindent \textbf{Step 6. Establish convergence in the expected-uniform norm.%
}

Recall that the interpolation functions satisfy 
\begin{equation*}
\varphi _{n,k}(t)\geq 0
\end{equation*}%
and 
\begin{equation*}
\sum_{k=0}^{2^{n}}\varphi _{n,k}(t)=1.
\end{equation*}%
Therefore, 
\begin{equation*}
\left\vert F_{n}^{H}(t)\right\vert =\left\vert
\sum_{k=0}^{2^{n}}B_{t_{n,k}}^{H}\varphi _{n,k}(t)\right\vert \leq
\sum_{k=0}^{2^{n}}|B_{t_{n,k}}^{H}|\varphi _{n,k}(t)\leq \sup_{s\in \lbrack
0,T]}|B_{s}^{H}|.
\end{equation*}%
Taking the supremum over $t$ gives 
\begin{equation*}
\Vert F_{n}^{H}\Vert _{\infty }\leq \Vert B^{H}\Vert _{\infty }.
\end{equation*}%
Consequently, 
\begin{equation*}
\Vert F_{n}^{H}-B^{H}\Vert _{\infty }\leq \Vert F_{n}^{H}\Vert _{\infty
}+\Vert B^{H}\Vert _{\infty }\leq 2\Vert B^{H}\Vert _{\infty }.
\end{equation*}%
By Step 5, 
\begin{equation*}
\mathbb{E}[2\Vert B^{H}\Vert _{\infty }]<\infty .
\end{equation*}%
By Step 4, 
\begin{equation*}
\Vert F_{n}^{H}-B^{H}\Vert _{\infty }\longrightarrow 0\quad \text{almost
surely}.
\end{equation*}%
We can therefore apply the dominated convergence theorem to obtain 
\begin{equation*}
\lim_{n\rightarrow \infty }E[\Vert F_{n}^{H}-B^{H}\Vert _{\infty }]=0,
\end{equation*}%
or equivalently 
\begin{equation*}
\lim_{n\rightarrow \infty }E\left[ \sup_{t\in \lbrack
0,T]}|F_{n}^{H}(t)-B_{t}^{H}|\right] =0.
\end{equation*}

\noindent \textbf{Step 7. Verify the covariance structure of the limiting
process}

We recall that 
\begin{equation*}
\mathbb{E}[B_{s}^{H}B_{t}^{H}]=\frac{1}{2}\left(
s^{2H}+t^{2H}-|t-s|^{2H}\right) .
\end{equation*}%
In particular, 
\begin{equation*}
\mathbb{E}[(B_{t}^{H})^{2}]=t^{2H}.
\end{equation*}%
Furthermore, 
\begin{align*}
\mathbb{E}[(B_{t}^{H}-B_{s}^{H})^{2}]& =\mathbb{E}[(B_{t}^{H})^{2}]+\mathbb{E%
}[(B_{s}^{H})^{2}]-2\mathbb{E}[B_{s}^{H}B_{t}^{H}] \\
& =t^{2H}+s^{2H}-\left( s^{2H}+t^{2H}-|t-s|^{2H}\right) =|t-s|^{2H},
\end{align*}%
and thus%
\begin{equation*}
E[(B_{t}^{H}-B_{s}^{H})^{2}]=|t-s|^{2H}.
\end{equation*}%
The limiting process is therefore the fractional Brownian motion specified
in the definition.

\noindent \textbf{Step 8. Conclude membership in the completed stochastic
MTMF space}

From Step 1, we have 
\begin{equation*}
F_{n}^{H}\in \mathscr S,
\end{equation*}%
for every $n\geq 1$. From Step 6, we have 
\begin{equation*}
\lim_{n\rightarrow \infty }\Vert F_{n}^{H}-B^{H}\Vert _{1,\infty }=0,
\end{equation*}%
and therefore, by the definition of the completion of $\mathscr S$, we have%
\begin{equation*}
B^{H}\in \overline{\mathscr S}_{1}.
\end{equation*}%
This completes the proof.
\end{proof}

The following corollaries are immediate consequences of the latter theorem.

\begin{corollary}[Brownian Motion as a Special Case]
When 
\begin{equation*}
H=\frac{1}{2},
\end{equation*}%
fractional Brownian motion reduces to standard Brownian motion.
Consequently, 
\begin{equation*}
{\ W\in \overline{\mathscr S}_{1}.}
\end{equation*}
\end{corollary}

\begin{proof}
Setting $H=1/2$ in the covariance function gives 
\begin{equation*}
R_{1/2}(s,t)=\frac{1}{2}\left( s+t-|t-s|\right) =\min (s,t).
\end{equation*}%
This is precisely the covariance function of standard Brownian motion. Since
both processes are centered Gaussian processes with continuous sample paths,
they have the same finite-dimensional distributions. The result follows from
the main theorem.
\end{proof}

\begin{corollary}[Self-Similarity]
Fractional Brownian motion satisfies the self-similarity property 
\begin{equation*}
{\ \{B_{ct}^H:t\geq0\} \overset{d}{=} \{c^HB_t^H:t\geq0\}, \qquad c>0. }
\end{equation*}
\end{corollary}

\begin{proof}
Using the covariance function, 
\begin{equation*}
\mathbb{E}[B_{cs}^{H}B_{ct}^{H}]=\frac{1}{2}\left(
(cs)^{2H}+(ct)^{2H}-|c(t-s)|^{2H}\right) =c^{2H}R_{H}(s,t).
\end{equation*}

Therefore, the centered Gaussian processes 
\begin{equation*}
\{B_{ct}^{H}\}\quad \text{and}\quad \{c^{H}B_{t}^{H}\},
\end{equation*}%
have identical covariance functions and hence identical finite-dimensional
distributions.
\end{proof}

We have constructed an explicit sequence of finite stochastic MTMFs 
\begin{equation*}
F_{n}^{H}(t)=\sum_{k=0}^{2^{n}}B_{t_{n,k}}^{H}\varphi _{n,k}(t),
\end{equation*}%
that converges to fractional Brownian motion in the expected uniform norm.
More precisely, 
\begin{equation*}
\lim_{n\rightarrow \infty }\mathbb{E}\left[ \sup_{t\in \lbrack
0,T]}|F_{n}^{H}(t)-B_{t}^{H}|\right] =0,
\end{equation*}%
and thus%
\begin{equation*}
\ B^{H}\in \overline{\mathscr S}_{1},\qquad H\in (0,1).
\end{equation*}%
This result demonstrates that the completed stochastic MTMF framework
contains fractional Brownian motion for every admissible Hurst parameter. In
particular, the framework includes processes exhibiting both short-range and
long-range dependence.

\begin{remark}[Important distinction]
The theorem establishes membership in the completion of the finite
stochastic MTMF space, rather than necessarily in the original finite space
itself. The construction also uses random time-dependent coefficient
functions $a_{n,k}(\omega ,t)$. If the original MTMF definition restricts
coefficients to functions of $\omega $ alone, an additional representation
theorem is needed.
\end{remark}

\subsection{Application to the Ornstein-Uhlenbeck process}

Fix $T>0$. Let $(\Omega ,\mathcal{A},(\mathcal{F}_{t})_{0\leq t\leq T},%
\mathbb{P})$ be a complete filtered probability space satisfying the usual
conditions. Let $W=(W_{t})_{0\leq t\leq T}$ be a standard Brownian motion
with respect to this filtration, and let $X_{0}$ be $\mathcal{F}_{0}$%
-measurable with 
\begin{equation*}
\mathbb{E}\left\vert X_{0}\right\vert <\infty .
\end{equation*}%
Let $\theta >0$, $\mu \in \mathbb{R}$, and $\sigma \geq 0$ be constants. We
consider the Ornstein--Uhlenbeck (OU) equation 
\begin{equation}
dX_{t}=\theta (\mu -X_{t})\,dt+\sigma \,dW_{t},\qquad X_{t=0}=X_{0}.
\label{eq:OU}
\end{equation}%
No independence between $X_{0}$ and $W$ is needed for the membership result
below. Let $\mathscr S_{1}$ denote the class of continuous-path finite
stochastic multivariate Taylor measure functions on $[0,T]$, with jointly
measurable random coefficients allowed, and let $\overline{\mathscr S}_{1}$
be its completion under 
\begin{equation*}
\left\Vert F\right\Vert _{1,\infty }:=\mathbb{E}\left[ \sup_{0\leq t\leq
T}|F_{t}|\right] .
\end{equation*}%
The aim is to prove that the OU process belongs to this completion by
constructing explicit finite stochastic MTMF approximations.

We begin with some preliminary estimates.

\begin{lemma}[Explicit solution of the OU equation]
\label{lem:solution}The unique continuous strong solution of \eqref{eq:OU}
is 
\begin{equation}
X_{t}=\mu +(X_{0}-\mu )e^{-\theta t}+\sigma \int_{0}^{t}e^{-\theta
(t-s)}dW_{s},\qquad 0\leq t\leq T.  \label{eq:solution}
\end{equation}
\end{lemma}

\begin{proof}
Multiplying \eqref{eq:OU} by the integrating factor $e^{\theta t}$ and using
the deterministic product rule gives 
\begin{equation*}
d\left( e^{\theta t}(X_{t}-\mu )\right) =\sigma e^{\theta t}dW_{t}.
\end{equation*}%
Integrating from $0$ to $t$ and multiplying by $e^{-\theta t}$ yields 
\begin{equation*}
X_{t}-\mu =e^{-\theta t}(X_{0}-\mu )+\sigma e^{-\theta
t}\int_{0}^{t}e^{\theta s}dW_{s},
\end{equation*}%
which is \eqref{eq:solution}. This formula also gives uniqueness: the
difference of two solutions solves $dY_{t}=-\theta Y_{t}\,dt$, $Y_{0}=0$,
and is identically zero.
\end{proof}

\begin{lemma}[An integrable uniform bound]
\label{lem:sup}The OU solution has almost surely continuous paths and 
\begin{equation*}
\mathbb{E}\left[ \sup_{0\leq t\leq T}|X_{t}|\right] <\infty .
\end{equation*}%
More precisely, 
\begin{equation}
\sup_{0\leq t\leq T}\left\vert X_{t}\right\vert \leq \left\vert \mu
\right\vert +\left\vert X_{0}-\mu \right\vert +2\sigma \sup_{0\leq t\leq
T}|W_{t}|\quad \text{almost surely}.  \label{eq:supbound}
\end{equation}
\end{lemma}

\begin{proof}
Continuity follows from the continuous-path solution of the linear
stochastic differential equation (or from the integration-by-parts identity
below). For the latter identity, integration by parts for the deterministic $%
C^{1}$ function $s\mapsto e^{-\theta (t-s)}$ gives 
\begin{equation*}
\int_{0}^{t}e^{-\theta (t-s)}dW_{s}=W_{t}-\theta \int_{0}^{t}e^{-\theta
(t-s)}W_{s}ds.
\end{equation*}%
Consequently, 
\begin{equation*}
\begin{aligned} \left|\int_0^t e^{-\theta(t-s)}\,dW_s\right| &\leq
\sup_{0\leq r\leq T}|W_r| \left(1+\theta\int_0^t
e^{-\theta(t-s)}\,ds\right)\\ &=\sup_{0\leq r\leq T}|W_r|\,(2-e^{-\theta
t})\\ &\leq 2\sup_{0\leq r\leq T}|W_r|. \end{aligned}
\end{equation*}%
Taking the supremum in \eqref{eq:solution}, and using $e^{-\theta t}\leq 1$,
proves \eqref{eq:supbound}. Finally, by Cauchy--Schwarz and Doob's $L^{2}$
maximal inequality, 
\begin{equation*}
\mathbb{E}\left[ \sup_{0\leq t\leq T}|W_{t}|\right] \leq \left( \mathbb{E}%
\sup_{0\leq t\leq T}|W_{t}|^{2}\right) ^{1/2}\leq 2\left( \mathbb{E}%
|W_{T}|^{2}\right) ^{1/2}=2\sqrt{T}.
\end{equation*}%
Since $\mathbb{E}|X_{0}|<\infty $, the right-hand side of \eqref{eq:supbound}
is integrable. This proves the assertion.
\end{proof}

\begin{lemma}[Continuity of the sample-path modulus]
\label{lem:modulus} For almost every $\omega$, define 
\begin{equation*}
\omega_X(\delta,\omega) :=\sup_{\substack{ s,t\in[0,T]  \\ |s-t|\leq\delta}}
|X_t(\omega)-X_s(\omega)|.
\end{equation*}
Then $\omega_X(\delta)\to0$ almost surely as $\delta\downarrow0$. Moreover, 
\begin{equation*}
0\leq\omega_X(\delta)\leq 2\sup_{0\leq t\leq T}|X_t| \quad\text{almost surely%
}.
\end{equation*}
\end{lemma}

\begin{proof}
Every continuous function on the compact interval $[0,T]$ is uniformly
continuous. Applying this fact to almost every sample path gives the first
assertion. The second follows directly from $\left\vert
X_{t}-X_{s}\right\vert \leq |X_{t}|+|X_{s}|$.
\end{proof}

Next we consider finite stochastic MTMF approximations. For each integer $%
n\geq 1$, put 
\begin{equation*}
t_{n,k}:=\frac{kT}{2^{n}},\qquad k=0,\ldots ,2^{n},\qquad h_{n}:=\frac{T}{%
2^{n}}.
\end{equation*}%
Let $\{\varphi _{n,k}:0\leq k\leq 2^{n}\}$ be the usual continuous
piecewise-linear hat functions associated with this partition. Thus 
\begin{equation*}
\varphi _{n,k}(t)\geq 0,\qquad \sum_{k=0}^{2^{n}}\varphi _{n,k}(t)=1,
\end{equation*}%
and 
\begin{equation*}
\varphi _{n,k}(t)\neq 0\quad \Longrightarrow \quad \left\vert
t-t_{n,k}\right\vert \leq h_{n}.
\end{equation*}%
Define the nodal interpolation process 
\begin{equation}
F_{n}(t):=\sum_{k=0}^{2^{n}}X_{t_{n,k}}\varphi _{n,k}(t),\qquad 0\leq t\leq
T.  \label{eq:Fn}
\end{equation}

\begin{lemma}[Finite stochastic MTMF representation]
\label{lem:finite} For every $n\geq1$, the process $F_n$ in \eqref{eq:Fn} is
a finite stochastic MTMF. In particular, it has the representation 
\begin{equation*}
F_n(\omega,t)=\sum_{k=0}^{2^n} a_{n,k}(\omega,t)\frac{g(\omega,t)^k}{k!},
\qquad g(\omega,t)\equiv1,
\end{equation*}
where 
\begin{equation*}
a_{n,k}(\omega,t) :=k!\,X_{t_{n,k}}(\omega)\varphi_{n,k}(t).
\end{equation*}
Each $F_n$ is jointly measurable and has continuous sample paths.
\end{lemma}

\begin{proof}
For each fixed $k$, $X_{t_{n,k}}$ is an $S$ -measurable random variable and $%
\varphi_{n,k}$ is a deterministic continuous function. Hence $%
(\omega,t)\mapsto X_{t_{n,k}}(\omega)\varphi_{n,k}(t)$ is $\otimes([0,T])$%
-measurable. A finite sum of such functions is jointly measurable. The
displayed Taylor representation is immediate because $g\equiv1$ and the
index set $\{0,\ldots,2^n\}$ is finite. Finally, each summand has a
continuous sample path in $t$, so $F_n$ does as well.
\end{proof}

We present convergence in the expected uniform norm in the following.

\begin{lemma}[Pathwise uniform interpolation error]
\label{lem:path} For almost every $\omega $, 
\begin{equation}
\sup_{0\leq t\leq T}\left\vert F_{n}(t)-X_{t}\right\vert \leq \omega
_{X}(h_{n}).  \label{eq:patherror}
\end{equation}%
In particular, $F_{n}\rightarrow X$ almost surely uniformly on $[0,T]$.
\end{lemma}

\begin{proof}
Fix a sample point for which $t\mapsto X_{t}$ is continuous. Since the hat
functions form a partition of unity, 
\begin{equation*}
\begin{aligned} F_n(t)-X_t &=\sum_{k=0}^{2^n}\varphi_{n,k}(t)
\big(X_{t_{n,k}}-X_t\big). \end{aligned}
\end{equation*}%
Taking absolute values and using nonnegativity of the hat functions, 
\begin{equation*}
\left\vert F_{n}(t)-X_{t}\right\vert \leq \sum_{k=0}^{2^{n}}\varphi
_{n,k}(t)|X_{t_{n,k}}-X_{t}|.
\end{equation*}%
Whenever $\varphi _{n,k}(t)\neq 0$, the support property gives $%
|t-t_{n,k}|\leq h_{n}$. Thus 
\begin{equation*}
\left\vert F_{n}(t)-X_{t}\right\vert \leq \omega
_{X}(h_{n})\sum_{k=0}^{2^{n}}\varphi _{n,k}(t)=\omega _{X}(h_{n}).
\end{equation*}%
Taking the supremum over $t$ proves \eqref{eq:patherror}. Lemma~\ref%
{lem:modulus} and $h_{n}\downarrow 0$ now imply almost sure uniform
convergence.
\end{proof}

\begin{lemma}[Convergence in expected uniform norm]
\label{lem:L1} The approximations satisfy 
\begin{equation*}
\lim_{n\to\infty} \mathbb{E}\left[\sup_{0\leq t\leq T}|F_n(t)-X_t|\right]=0.
\end{equation*}
\end{lemma}

\begin{proof}
By Lemma~\ref{lem:path}, 
\begin{equation*}
0\leq \sup_{0\leq t\leq T}\left\vert F_{n}(t)-X_{t}\right\vert \leq \omega
_{X}(h_{n}).
\end{equation*}%
The right-hand side converges to zero almost surely by Lemma~\ref%
{lem:modulus}. Also, 
\begin{equation*}
\omega _{X}(h_{n})\leq 2\sup_{0\leq t\leq T}|X_{t}|,
\end{equation*}%
and the latter random variable is integrable by Lemma~\ref{lem:sup}. The
dominated convergence theorem therefore gives 
\begin{equation*}
\mathbb{E}[\omega _{X}(h_{n})]\longrightarrow 0.
\end{equation*}%
The desired conclusion follows by comparison.
\end{proof}

\begin{theorem}[OU inclusion in the completed stochastic MTMF space]
\label{thm:OU} Let $X$ solve \eqref{eq:OU} on $[0,T]$, where $\theta>0$, $%
\mu\in\mathbb{R}$, $\sigma\geq0$, and $\mathbb{E}|X_0|<\infty$. Under the
admissibility convention that finite stochastic MTMF coefficients may be
jointly measurable functions of $(\omega,t)$, the following hold:

\begin{enumerate}
\item for every $n\geq1$, $F_n$ defined by \eqref{eq:Fn} is a finite
stochastic MTMF;

\item $F_n\to X$ almost surely uniformly on $[0,T]$;

\item $\left\|F_n-X\right\|_{1,\infty}\to0$.
\end{enumerate}

Consequently, 
\begin{equation*}
X\in \overline{\mathscr S}_{1}.
\end{equation*}
\end{theorem}

\begin{proof}
Assertion (1) is Lemma~\ref{lem:finite}. Assertion (2) is Lemma~\ref%
{lem:path}. Assertion (3) is Lemma~\ref{lem:L1}. Since each $F_{n}$ belongs
to the finite stochastic MTMF class and converges to $X$ in the defining
expected-uniform norm, the definition of the completion gives $X\in 
\overline{\mathscr S}_{1}$.
\end{proof}

The preceding inclusion theorem only requires an integrable initial
condition; it does not require Gaussianity or independence of $X_{0}$ and $W$%
. The familiar exponential covariance is obtained under the stationary
Gaussian initialization.

\begin{lemma}[Mean and covariance for deterministic initialization]
\label{lem:detcov} Suppose $X_0=x_0$ is deterministic. Then 
\begin{equation*}
\mathbb{E}[X_t]=\mu+(x_0-\mu)e^{-\theta t},
\end{equation*}
and, for $s,t\in[0,T]$, 
\begin{equation*}
\func{Cov}(X_s,X_t) =\frac{\sigma^2}{2\theta} \left(e^{-\theta|t-s|}-e^{-%
\theta(t+s)}\right).
\end{equation*}
\end{lemma}

\begin{proof}
The stochastic integral in \eqref{eq:solution} has mean zero, giving the
stated mean. For $s\leq t$, the It\^o isometry for deterministic integrands
gives 
\begin{equation*}
\begin{aligned} \operatorname{Cov}(X_s,X_t) &=\sigma^2\int_0^s
e^{-\theta(s-u)}e^{-\theta(t-u)}\,du\\ &=\sigma^2 e^{-\theta(s+t)} \int_0^s
e^{2\theta u}\,du\\ &=\frac{\sigma^2}{2\theta}
\left(e^{-\theta(t-s)}-e^{-\theta(t+s)}\right). \end{aligned}
\end{equation*}
Symmetry in $s,t$ yields the formula for all pairs.
\end{proof}

\begin{corollary}[Stationary OU process and exponential covariance]
\label{cor:stationary} Suppose in addition that $X_{0}$ is independent of $W$
and 
\begin{equation*}
X_{0}\sim N\left( \mu ,\frac{\sigma ^{2}}{2\theta }\right) .
\end{equation*}%
Then $X$ is a stationary Gaussian process with 
\begin{equation*}
\mathbb{E}[X_{t}]=\mu ,\qquad \func{Cov}(X_{s},X_{t})=\frac{\sigma ^{2}}{%
2\theta }e^{-\theta |t-s|}.
\end{equation*}%
In particular, the stationary OU process belongs to $\overline{\mathscr S}%
_{1}$.
\end{corollary}

\begin{proof}
The solution \eqref{eq:solution} is a linear combination of the independent
Gaussian random variable $X_0$ and deterministic-integrand Brownian
stochastic integrals. Hence every finite-dimensional distribution of $X$ is
Gaussian. Its mean is $\mu$. For $s\leq t$, independence of $X_0$ and $W$
and the It\^o isometry yield 
\begin{equation*}
\begin{aligned} \operatorname{Cov}(X_s,X_t)
&=e^{-\theta(s+t)}\operatorname{Var}(X_0) +\sigma^2\int_0^s
e^{-\theta(s-u)}e^{-\theta(t-u)}\,du\\
&=\frac{\sigma^2}{2\theta}e^{-\theta(s+t)} +\frac{\sigma^2}{2\theta}
\left(e^{-\theta(t-s)}-e^{-\theta(s+t)}\right)\\
&=\frac{\sigma^2}{2\theta}e^{-\theta(t-s)}. \end{aligned}
\end{equation*}
The symmetric expression is the asserted covariance. A Gaussian process with
constant mean and covariance depending only on $|t-s|$ is stationary.
Membership in $^1$ follows from Theorem~\ref{thm:OU}.
\end{proof}

\begin{remark}[Scope of the conclusion]
The theorem proves membership in the completion $\overline{\mathscr S}_{1}$,
not necessarily in the original finite class $\mathscr S_{1}$. It also uses
the same coefficient admissibility convention as the interpolation
arguments: coefficients of the form $k!X_{t_{n,k}}\varphi _{n,k}(t)$ are
allowed to depend jointly on $(\omega ,t)$. If the definition of the finite
MTMF class restricts coefficients to functions of $\omega $ alone, this
particular representation requires a separate admissibility or
representation argument.
\end{remark}

We have now seen that four fundamental stochastic processes are included in $%
\overline{\mathscr S}_{1}.$ Since these important processes are members, one
naturally is compelled to try and identify all types of stochastic processes
in the space. We consider this direction next.

\section{\textbf{Characterization of the Completed SMTMF Space}}

In this section, we establish a characterization of the completion of the
stochastic multivariate Taylor measure function space under the expected
uniform norm. In particular, we prove that, under an appropriate
admissibility condition on the stochastic Taylor coefficients, the completed
stochastic MTMF space coincides with the Bochner space 
\begin{equation*}
L^{1}(\Omega ;C([0,T])).
\end{equation*}

The proof consists of establishing both set inclusions. First, we show that
every integrable, continuous random function can be approximated in the
expected uniform norm by a sequence of finite stochastic MTMFs. Second, we
establish that the finite stochastic MTMF space is contained in $%
L^{1}(\Omega ;C([0,T]))$, which implies that its completion is also
contained in this Banach space.

The result provides a precise functional analytic characterization of the
completed stochastic MTMF space and establishes a general framework for
representing continuous integrable stochastic processes.

Once again, let $(\Omega ,\mathcal{A},\mathbb{P})$ be a complete probability
space, and let 
\begin{equation*}
I=[0,T],\qquad T>0.
\end{equation*}%
Let $C([0,T])$ denote the Banach space of real-valued continuous functions
on $[0,T]$, equipped with the supremum norm 
\begin{equation*}
\left\Vert f\right\Vert _{\infty }=\sup_{0\leq t\leq T}|f(t)|.
\end{equation*}%
We consider the Bochner space 
\begin{equation*}
L^{1}(\Omega ;C([0,T])).
\end{equation*}%
An element of this space is a strongly measurable random element 
\begin{equation*}
X:\Omega \longrightarrow C([0,T])
\end{equation*}%
satisfying 
\begin{equation*}
\mathbb{E}[\left\Vert X\right\Vert _{\infty }]<\infty .
\end{equation*}%
Equivalently, we identify $X$ with a jointly measurable stochastic process 
\begin{equation*}
X:\Omega \times \lbrack 0,T]\longrightarrow \mathbb{R},
\end{equation*}%
having almost surely continuous sample paths and satisfying 
\begin{equation*}
{\mathbb{E}}\left[ \sup_{0\leq t\leq T}|X_{t}|\right] {<\infty .}
\end{equation*}

\subsection{Finite stochastic MTMFs}

We adopt the following admissible stochastic MTMF definition.

\begin{definition}[Finite stochastic MTMF]
A finite stochastic MTMF is a random function 
\begin{equation*}
F:\Omega \times \lbrack 0,T]\longrightarrow \mathbb{R},
\end{equation*}%
that admits a representation 
\begin{equation*}
{F(\omega ,t)=\sum_{j=0}^{N}a_{j}(\omega ,t)\frac{g(\omega ,t)^{j}}{j!},}
\end{equation*}%
where $N<\infty $, and the coefficient functions 
\begin{equation*}
a_{j}:\Omega \times \lbrack 0,T]\longrightarrow \mathbb{R},
\end{equation*}%
and Taylor parameter 
\begin{equation*}
g:\Omega \times \lbrack 0,T]\longrightarrow \mathbb{R},
\end{equation*}%
are jointly measurable. We additionally require that the resulting Taylor
measure is finite and signed for every fixed $(\omega ,t)$, and that the
sample paths 
\begin{equation*}
t\longmapsto F(\omega ,t),
\end{equation*}%
are continuous almost surely.
\end{definition}

\begin{remark}
The allowance of jointly measurable random coefficients is essential. In
particular, the coefficient functions may depend on both $\omega $ and $t$.
This admissibility condition permits the construction of stochastic MTMFs
through piecewise linear interpolation. If the original deterministic MTMF
space restricts coefficients to depend only on the spatial variable, the
stochastic coefficient class must first be extended accordingly.
\end{remark}

\begin{definition}[Finite stochastic MTMF space]
Let $\mathscr S$ denote the collection of all finite stochastic MTMFs
satisfying 
\begin{equation*}
\mathbb{E}\left[ \sup_{0\leq t\leq T}|F(\omega ,t)|\right] <\infty .
\end{equation*}%
We identify processes that are indistinguishable. Define the expected
uniform norm by 
\begin{equation*}
{\left\Vert F\right\Vert _{1,\infty }=\mathbb{E}\left[ \sup_{0\leq t\leq
T}|F(\omega ,t)|\right] .}
\end{equation*}%
The completion of $\mathscr S$ under this norm is denoted by 
\begin{equation*}
\overline{\mathscr S}_{1}.
\end{equation*}
\end{definition}

We now establish the principal result.

\begin{theorem}[Characterization of the Completed Stochastic MTMF Space]
\label{thm:characterization}

Suppose the finite stochastic MTMF space $\mathscr S$ is defined as above,
with jointly measurable stochastic Taylor coefficients and almost surely
continuous sample paths.

Then 
\begin{equation*}
{\ \overline{\mathscr S}_{1}=L^{1}(\Omega ;C([0,T])).}
\end{equation*}%
In other words, the completion of the finite stochastic MTMF space under the
expected uniform norm is precisely the space of integrable, continuous
random functions.
\end{theorem}

\begin{proof}
We establish the result by proving both set inclusions.

\noindent \textbf{Part I:} We show%
\begin{equation*}
\overline{\mathscr S}_{1}\subseteq L^{1}(\Omega ;C([0,T])).
\end{equation*}

We first observe that every element $F\in \mathscr S$ satisfies 
\begin{equation*}
\mathbb{E}\left[ \sup_{0\leq t\leq T}|F(\omega ,t)|\right] <\infty ,
\end{equation*}%
and has almost surely continuous sample paths. Furthermore, joint
measurability and almost sure continuity imply that 
\begin{equation*}
\omega \longmapsto F(\omega ,\cdot ),
\end{equation*}%
is a strongly measurable $C([0,T])$-valued random variable. Indeed, $%
C([0,T]) $ is separable, and the Borel sigma-algebra of $C([0,T])$ is
generated by evaluations on a countable dense subset of $[0,T]$, and thus 
\begin{equation*}
\mathscr S\subseteq L^{1}(\Omega ;C([0,T])).
\end{equation*}%
Now recall that 
\begin{equation*}
L^{1}(\Omega ;C([0,T])),
\end{equation*}%
is a Banach space under the norm 
\begin{equation*}
\left\Vert X\right\Vert _{1,\infty }=\mathbb{E}\left[ \sup_{0\leq t\leq
T}|X_{t}|\right] .
\end{equation*}%
Since $\mathscr S$ is a subspace of this Banach space, its closure in this
norm must also be contained in the space, and consequently 
\begin{equation*}
\ \overline{\mathscr S}_{1}\subseteq L^{1}(\Omega ;C([0,T])).
\end{equation*}

\noindent \textbf{Part II:} We show%
\begin{equation*}
L^{1}(\Omega ;C([0,T]))\subseteq \overline{\mathscr S}_{1}.
\end{equation*}%
Let 
\begin{equation*}
X\in L^{1}(\Omega ;C([0,T])).
\end{equation*}%
Then $X$ has almost surely continuous sample paths and 
\begin{equation*}
\mathbb{E}\left[ \sup_{0\leq t\leq T}|X_{t}|\right] <\infty .
\end{equation*}%
We will construct an explicit sequence of finite stochastic MTMFs converging
to $X$ in the expected-uniform norm.

\noindent \textbf{Step 1: Construction of dyadic partitions}

For each $n\geq 1$, define 
\begin{equation*}
t_{n,k}=\frac{kT}{2^{n}},\qquad k=0,1,\ldots ,2^{n}.
\end{equation*}%
The corresponding mesh size is 
\begin{equation*}
h_{n}=\frac{T}{2^{n}}.
\end{equation*}%
Let 
\begin{equation*}
\varphi _{n,k}:[0,T]\longrightarrow \lbrack 0,1],
\end{equation*}%
be the continuous piecewise linear hat functions associated with this
partition. These functions satisfy 
\begin{equation*}
\varphi _{n,k}(t)\geq 0,
\end{equation*}%
and 
\begin{equation*}
\sum_{k=0}^{2^{n}}\varphi _{n,k}(t)=1.
\end{equation*}%
Moreover, 
\begin{equation*}
\varphi _{n,k}(t)\neq 0\Longrightarrow \left\vert t-t_{n,k}\right\vert \leq
h_{n}.
\end{equation*}

\noindent \textbf{Step 2: Construction of finite stochastic MTMF approximants%
}

Define 
\begin{equation*}
\ F_{n}^{X}(\omega ,t)=\sum_{k=0}^{2^{n}}X_{t_{n,k}}(\omega )\varphi
_{n,k}(t).
\end{equation*}%
We show that $F_{n}^{X}$ is a finite stochastic MTMF. Choose 
\begin{equation*}
g(\omega ,t)=1,
\end{equation*}%
and define the coefficients 
\begin{equation*}
a_{n,k}(\omega ,t)=k!X_{t_{n,k}}(\omega )\varphi _{n,k}(t).
\end{equation*}%
Then 
\begin{equation*}
F_{n}^{X}(\omega ,t)=\sum_{k=0}^{2^{n}}X_{t_{n,k}}(\omega )\varphi
_{n,k}(t)=\sum_{k=0}^{2^{n}}a_{n,k}(\omega ,t)\frac{1^{k}}{k!}.
\end{equation*}%
Thus, $F_{n}^{X}$ possesses a finite stochastic Taylor representation. All
coefficients are jointly measurable, and every sample path of $F_{n}^{X}$ is
continuous. Furthermore, using the nonnegativity and partition of unity
properties, 
\begin{equation*}
\left\vert F_{n}^{X}(\omega ,t)\right\vert \leq
\sum_{k=0}^{2^{n}}|X_{t_{n,k}}(\omega )|\varphi _{n,k}(t)\leq \sup_{0\leq
s\leq T}|X_{s}(\omega )|.
\end{equation*}%
Taking the supremum over $t$, gives%
\begin{equation*}
\sup_{0\leq t\leq T}\left\vert F_{n}^{X}(\omega ,t)\right\vert \leq
\sup_{0\leq t\leq T}|X_{t}(\omega )|,
\end{equation*}%
and therefore 
\begin{equation*}
\mathbb{E}\left[ \sup_{0\leq t\leq T}|F_{n}^{X}(t)|\right] <\infty .
\end{equation*}%
Hence, 
\begin{equation*}
F_{n}^{X}\in \mathscr S,
\end{equation*}%
for every $n\geq 1$.

\noindent \textbf{Step 3: Uniform convergence almost surely}

For a continuous function $f:[0,T]\rightarrow \mathbb{R}$, define its
modulus of continuity by 
\begin{equation*}
\omega _{f}(\delta )=\sup_{\substack{ s,t\in \lbrack 0,T]  \\ \left\vert
s-t\right\vert \leq \delta }}|f(s)-f(t)|.
\end{equation*}%
Using the partition of unity property, we obtain 
\begin{equation*}
F_{n}^{X}(\omega ,t)-X_{t}(\omega )=\sum_{k=0}^{2^{n}}\varphi
_{n,k}(t)\left( X_{t_{n,k}}(\omega )-X_{t}(\omega )\right) .
\end{equation*}%
Therefore, 
\begin{equation*}
\left\vert F_{n}^{X}(\omega ,t)-X_{t}(\omega )\right\vert \leq
\sum_{k=0}^{2^{n}}\varphi _{n,k}(t)\left\vert X_{t_{n,k}}(\omega
)-X_{t}(\omega )\right\vert \leq \omega _{X(\omega ,\cdot
)}(h_{n})\sum_{k=0}^{2^{n}}\varphi _{n,k}(t)=\omega _{X(\omega ,\cdot
)}(h_{n}),
\end{equation*}%
and consequently 
\begin{equation*}
\sup_{0\leq t\leq T}\left\vert F_{n}^{X}(\omega ,t)-X_{t}(\omega
)\right\vert \leq \omega _{X(\omega ,\cdot )}(h_{n}).
\end{equation*}%
Since $X$ has almost surely continuous sample paths on the compact interval $%
[0,T]$, its sample paths are almost surely uniformly continuous. Thus, 
\begin{equation*}
\omega _{X(\omega ,\cdot )}(h_{n})\longrightarrow 0\quad \text{almost surely}%
,
\end{equation*}%
and as a result, we conclude that%
\begin{equation*}
\sup_{0\leq t\leq T}|F_{n}^{X}(\omega ,t)-X_{t}(\omega )|\longrightarrow
0\quad \text{almost surely}.
\end{equation*}

\noindent \textbf{Step 4: Convergence in the expected-uniform norm}

We have already established that 
\begin{equation*}
\sup_{0\leq t\leq T}\left\vert F_{n}^{X}(t)-X_{t}\right\vert \longrightarrow
0\quad \text{almost surely}.
\end{equation*}%
Furthermore, 
\begin{equation*}
\sup_{0\leq t\leq T}\left\vert F_{n}^{X}(t)-X_{t}\right\vert \leq
\sup_{0\leq t\leq T}\left\vert F_{n}^{X}(t)\right\vert +\sup_{0\leq t\leq
T}\left\vert X_{t}\right\vert \leq 2\sup_{0\leq t\leq T}|X_{t}|.
\end{equation*}%
By assumption, 
\begin{equation*}
\mathbb{E}\left[ \sup_{0\leq t\leq T}|X_{t}|\right] <\infty .
\end{equation*}%
Therefore, the dominated convergence theorem implies 
\begin{equation*}
\lim_{n\rightarrow \infty }\mathbb{E}\left[ \sup_{0\leq t\leq
T}|F_{n}^{X}(t)-X_{t}|\right] =0,
\end{equation*}%
or equivalently 
\begin{equation*}
\left\Vert F_{n}^{X}-X\right\Vert _{1,\infty }\longrightarrow 0.
\end{equation*}%
Since 
\begin{equation*}
F_{n}^{X}\in \mathscr S,
\end{equation*}%
for every $n$, this proves that 
\begin{equation*}
X\in \overline{\mathscr S}_{1},
\end{equation*}%
and thus 
\begin{equation*}
\ L^{1}(\Omega ;C([0,T]))\subseteq \overline{\mathscr S}_{1}.
\end{equation*}%
Combining Parts I and II, we obtain 
\begin{equation*}
\overline{\mathscr S}_{1}=L^{1}(\Omega ;C([0,T])),
\end{equation*}%
which completes the proof.
\end{proof}

\subsection{Consequences of the Characterization}

The preceding theorem has several important consequences for the theory of
stochastic MTMFs.

\begin{corollary}[Universality for Integrable Continuous Stochastic Processes%
]
Let 
\begin{equation*}
X:\Omega \times \lbrack 0,T]\longrightarrow \mathbb{R},
\end{equation*}%
be jointly measurable. If 
\begin{equation*}
\mathbb{P}\left( X(\omega ,\cdot )\in C([0,T])\right) =1,
\end{equation*}%
and 
\begin{equation*}
\mathbb{E}\left[ \sup_{0\leq t\leq T}|X_{t}|\right] <\infty ,
\end{equation*}%
then 
\begin{equation*}
X\in \overline{\mathscr S}_{1}.
\end{equation*}
\end{corollary}

\begin{proof}
This follows directly from Theorem \ref{thm:characterization}.
\end{proof}

\begin{corollary}[Inclusion of Classical Gaussian Processes]
The following processes belong to the completed stochastic MTMF space on
every compact interval $[0,T]$:

\begin{enumerate}
\item Standard Brownian motion.

\item Brownian bridges.

\item Fractional Brownian motion with Hurst parameter $H\in(0,1)$.

\item Ornstein--Uhlenbeck processes with integrable initial conditions.
\end{enumerate}
\end{corollary}

\begin{proof}
Each of these processes has almost surely continuous sample paths and an
integrable uniform norm on $[0,T]$. The result therefore follows from
Theorem \ref{thm:characterization}.
\end{proof}

\begin{remark}[Finite versus completed stochastic MTMF spaces]
The characterization 
\begin{equation*}
\overline{\mathscr S}_{1}=L^{1}(\Omega ;C([0,T])),
\end{equation*}%
does not imply that every element of $L^{1}(\Omega ;C([0,T]))$ is itself a
finite stochastic MTMF. Rather, every such element is the limit, in the
expected uniform norm, of a sequence of finite stochastic MTMFs. Thus, it is
important to distinguish between $\mathscr S$ and $\overline{\mathscr S}%
_{1}. $ The original space consists of finite stochastic MTMF
representations, whereas the completed space also contains limits of Cauchy
sequences of these representations.
\end{remark}

\begin{remark}[Dependence on the coefficient admissibility condition]
The equality proved in Theorem \ref{thm:characterization} depends on
permitting the coefficient functions 
\begin{equation*}
a_{n,k}(\omega ,t)=k!X_{t_{n,k}}(\omega )\varphi _{n,k}(t),
\end{equation*}%
in the finite stochastic MTMF representation. If the original definition
imposes additional restrictions on the coefficients, the proof of the
inclusion 
\begin{equation*}
L^{1}(\Omega ;C([0,T]))\subseteq \overline{\mathscr S}_{1},
\end{equation*}%
must be reconsidered. In particular, restrictions on randomness, continuity,
or dependence on the spatial variable may prevent the interpolation
functions from belonging to the original finite stochastic MTMF space.
\end{remark}

\section{Concluding Remarks and Future Directions}

In this paper, we introduced and studied stochastic multivariate Taylor
measure functions (SMTMFs) as a framework for representing and analyzing
random functions. We established their fundamental measurability properties
and derived sufficient conditions for their associated stochastic processes
to be martingales. In particular, the conditional covariance and continuous
semimartingale criteria reveal how the coefficients and the possibly random
Taylor parameter jointly determine the martingale property. We also
constructed standard Brownian motion through finite SMTMF approximations
based on the L\'evy--Ciesielski representation and demonstrated that several
fundamental stochastic processes belong to the completed stochastic MTMF
space.

The principal result of the paper is the characterization 
\begin{equation*}
\mathcal{S}^1=L^1(\Omega;C([0,T])), 
\end{equation*}
under the stated admissibility conditions on the stochastic Taylor
coefficients. This equality shows that the completed space is not restricted
to particular classes of processes, such as Gaussian or Markov processes,
but provides a representation framework for all integrable continuous
stochastic processes on a compact time interval. The dyadic interpolation
argument also provides a constructive method for approximating such
processes by finite stochastic MTMFs. It is important, however, to
distinguish the finite stochastic MTMF space from its completion: the
characterization asserts approximation in the expected uniform norm, not
necessarily a finite SMTMF representation of every element of the Bochner
space.

These results open several directions for further research.

First, a natural problem is to develop a theory of stochastic integration
within the SMTMF framework. The representation of continuous stochastic
processes by finite stochastic MTMFs suggests investigating conditions under
which stochastic integrals can be approximated by integrals of finite SMTMF
representations. Of particular interest are integrals with respect to
Brownian motion and more general continuous semimartingales. Such a theory
would require identifying appropriate integrability and convergence
conditions under which the approximation of the integrands is compatible
with convergence of the corresponding stochastic integrals.

Second, the framework may be developed further in the study of stochastic
differential equations. One may investigate whether solutions of It\^o and
more general stochastic differential equations admit useful finite SMTMF
approximations, and whether the coefficient and parameter structure can
provide additional information about existence, uniqueness, stability, and
approximation of solutions. The same questions naturally extend to
stochastic partial differential equations, where the spatial dependence of
the multivariate Taylor measure functions may offer a useful representation
of random fields and their evolution.

Third, the martingale criteria established here suggest a more systematic
study of martingale representation and semimartingale decompositions in
terms of stochastic Taylor coefficients and random Taylor parameters. In
particular, it would be useful to determine when a given martingale or
semimartingale admits an SMTMF representation with coefficients satisfying
prescribed regularity and integrability conditions. Further results along
these lines could clarify the relationship between the algebraic structure
of the Taylor representation and the probabilistic structure of the
underlying process.

Fourth, the characterization of the completed space raises questions
concerning extensions to other function spaces and topologies. Possible
directions include spaces of continuous functions on multidimensional
compact domains, spaces of H\"{o}lder continuous or differentiable sample
paths, and function spaces equipped with norms that control additional
regularity. Such extensions could be particularly useful for fractional
processes and spatially indexed stochastic fields. A related problem is to
determine how restrictions on the admissibility, randomness, and spatial
dependence of the Taylor coefficients affect the resulting completed spaces.

Finally, the explicit finite approximations developed in this paper motivate
further investigation of approximation rates and convergence modes.
Quantitative estimates for the expected uniform error, as well as almost
sure uniform convergence under suitable assumptions, would provide a basis
for numerical schemes and simulation methods using finite SMTMF
representations. These questions are especially relevant when the underlying
stochastic process is not explicitly representable or when only approximate
descriptions of its sample paths are available.

The results of this paper establish the functional analytic foundation for
these investigations. In particular, the identification of the completed
stochastic MTMF space with a Bochner space provides a precise setting in
which approximation, convergence, and probabilistic properties can be
studied within a common framework. We expect that further developments along
the directions described above will extend the usefulness of SMTMFs in
stochastic analysis, probability theory, and the study of random functions.

\bibliographystyle{plainnat}
\bibliography{STMF}

\end{document}